\documentclass[10pt,reqno]{amsart}

\usepackage{amsmath,amsfonts}
\usepackage{subfigure}
\usepackage{graphicx}
\usepackage{color}
\usepackage{tikz}
\usepackage{url}
\usepackage{xr}
\usepackage{etoolbox}
\AtBeginEnvironment{algorithm}{\let\textbf\relax}
\usepackage{enumerate}
\usepackage[pagebackref,colorlinks,citecolor=blue,linkcolor=magenta]{hyperref}

\usepackage[ruled,vlined]{algorithm2e}

\newtheorem{theorem}{Theorem}

\newtheorem{lemma}[theorem]{Lemma}
\newtheorem{proposition}[theorem]{Proposition}

\theoremstyle{definition}

\newtheorem{assump}[theorem]{Assumption}

 \theoremstyle{remark}
 
 \newtheorem{example}[theorem]{Example}

 \numberwithin{equation}{section}

\DeclareMathOperator{\CI}{CI}

\DeclareMathOperator{\argmax}{argmax}

\DeclareMathOperator{\score}{score}
\DeclareMathOperator{\parents}{Pa}
\DeclareMathOperator{\descendants}{De}
\DeclareMathOperator{\nondescendants}{Nd}
\DeclareMathOperator{\ancestors}{An}

\DeclareMathOperator*{\pa}{Pa}
\DeclareMathOperator*{\adj}{adj}

\newcommand{\R}{\mathbb{R}}
\newcommand{\CC}{\mathcal{C}}
\newcommand{\p}{\mathbb{P}}

\newcommand{\G}{\mathcal{G}}
\newcommand{\HH}{\mathcal{H}}
\newcommand{\M}{\mathcal{M}}

\newcommand\independent{\protect\mathpalette{\protect\independenT}{\perp}}
\def\independenT#1#2{\mathrel{\rlap{$#1#2$}\mkern2mu{#1#2}}}

\newcommand{\union}{\cup}
\newcommand{\intersection}{\cap}
\newcommand{\PP}{\mathbb{P}}

\SetKwInput{KwInput}{Input}
\SetKwInput{KwOutput}{Output}

\keywords{causal inference; Bayesian network; DAG; Gaussian graphical model; greedy search; pointwise and high-dimensional consistency; faithfulness; generalized permutohedron; DAG associahedron.}

\begin{document}

\title[Ordering-based Causal Inference]{Consistency Guarantees for Greedy Permutation-Based Causal Inference Algorithms}

\date{\today}

\author{Liam Solus}
\address{Institutionen f\"or Matematik, KTH, SE-100 44 Stockholm, Sweden}
\email{solus@kth.se}

\author{Yuhao Wang}
\address{Institute for Interdisciplinary Information Sciences, Tsinghua University, Beijing, China, 
and Shanghai Qi Zhi Institute,
Shanghai, China.}
\email{yuhaow@tsinghua.edu.cn}

\author{Caroline Uhler}
\address{Department of Electrical Engineering and Computer Science, and 
	Institute for Data, Systems and Society, Massachusetts Institute of Technology, Cambridge, MA, USA.}
\email{cuhler@mit.edu}

\begin{abstract}
Directed acyclic graphical models, or DAG models,  are widely used to represent complex causal systems.  
Since the basic task of learning such a model from data is NP-hard, a standard approach is greedy search over the space of directed acyclic graphs or Markov equivalence classes of directed acyclic graphs. 
As the space of directed acyclic graphs on $p$ nodes and the associated space of Markov equivalence classes are both much larger than the space of permutations, it is desirable to consider permutation-based greedy searches. 
Here, we provide the first consistency guarantees, both uniform and high-dimensional, of a greedy permutation-based search. 
This search corresponds to a simplex-like algorithm operating over the edge-graph of a sub-polytope of the permutohedron, called a DAG associahedron. 
Every vertex in this polytope is associated with a directed acyclic graph, and hence with a collection of permutations that are consistent with the directed acyclic graph ordering. 
A walk is performed on the edges of the polytope maximizing the sparsity of the associated directed acyclic graphs. 
We show via simulated and real data that this permutation search is competitive with current approaches.
\end{abstract}

\maketitle
\thispagestyle{empty}

\section{Introduction}
\label{sec: introduction}
Bayesian networks, or DAG models, are widely used to model complex causal systems arising, for example, in computational biology, epidemiology, or sociology \cite{Friedman_2000,Pearl_2000,Robins_2000,Spirtes_2001}.  
Given a directed acyclic graph, i.e., a DAG, $\G := ([p], A)$ with node set $[p] := \{1,2,\ldots,p\}$ and arrow set $A$, a DAG model associates to each node $i\in[p]$ of $\G$ a random variable $X_i$.  
By the Markov property, $\G$ encodes a set of conditional independence relations
$
X_i \independent X_{\nondescendants(i)\backslash\parents(i)} \,\vert\, X_{\parents(i)},
$
where $\nondescendants(i)$ and $\parents(i)$, respectively, denote the nondesendants and parents of the node $i$ in $\G$.  
A joint distribution $\mathbb{P}$ on $X_1,\ldots,X_p$ is said to satisfy the Markov assumption, or be Markov, with respect to $\mathcal{G}$ if it entails these conditional independence relations. 
This paper is concerned with the structural learning problem:  Suppose we sample from a distribution $\p$ that is Markov with respect to a DAG $\G^\ast$.  If we infer from this data a collection of conditional independence relations $\CC$, can we recover the unknown DAG $\G^\ast$ using $\CC$? 

In general, this problem is not well-defined since multiple DAGs can encode the same set of conditional independence relations.  
Any two such DAGs are termed Markov equivalent, and they are said to belong to the same Markov equivalence class.  
Thus, our goal becomes to identify the Markov equivalence class $\M(\G^\ast)$ of $\G^\ast$. 
The Markov assumption alone is not sufficient to guarantee identifiability, and so additional identifiability assumptions have been studied, the most prominent being faithfulness \cite{Spirtes_2001}.
Unfortunately, this assumption has been shown to be restrictive in practice~\cite{URBY13}. 
Hence, it is desirable to develop structure learning algorithms that are consistent under strictly weaker assumptions than faithfulness. 

Since the space of all Markov equivalence classes of DAGs on $p$ nodes grows super-exponentially in $p$ \cite{GP01}, one way to perform structure learning is to use greedy approaches. 
For example, the greedy equivalence search \cite{C02,M97} greedily maximizes a score, such as the Bayesian information criterion, over the space of all Markov equivalence classes on $p$ nodes.  
An alternative approach is to consider algorithms with a reduced search space, such as the space of all $p!$ linear extensions of DAGs; i.e., the permutations of $[p]$.  
Greedy permutation-based algorithms combine both of these heuristic approaches for DAG model learning.
In recent decades, a variety of greedy permutation-based algorithms have been proposed and analyzed. 
See, for instance, \cite{B92,CH92,L96,SV93,TK12}.  
However, all of these algorithms rely on heuristics for sparse DAG recovery and are therefore not provably consistent, even under the faithfulness assumption. 

On the other hand, a non-greedy permutation-based algorithm known as the sparsest permutation algorithm was introduced in \cite{RU13}, and it was shown to be consistent under strictly weaker assumptions than faithfulness.  
Unfortunately, the sparsest permutation algorithm must generate and score a DAG for each permutation of $\{1,\ldots,p\}$, and hence it runs in $\mathcal{O}(p!)$ time no matter the true underlying DAG model.
Here, we provide the first consistency guarantees of a greedy permutation-based algorithm for DAG model structure learning. 
This algorithm is a greedy version of the sparsest permutation algorithm. 
Unlike its non-greedy predecessor, the proposed algorithm scales to DAG structure discovery with hundreds of variables, since only in the worst case does it have to search over all $p!$ permutations; e.g., when the true model is the complete graph. 
Such worst-case behavior is to be expected since in general the problem of learning a DAG model is NP-hard~\cite{Chickering_NPhard}.
In addition, we show that our greedy sparsest permutation algorithm is consistent under weaker assumptions than standardly assumed, which translates into competitive performance in terms of structure recovery when compared to currently popular algorithms on both simulated and real data.

\section{Background}
\label{sec: identifiability assumptions}
We refer the reader to Section~\ref{app: greedy SP} in the Supplementary Material for graph theory definitions and notation.  
A fundamental result about DAG models is that the complete set of conditional independence relations implied by the Markov assumption for $\G$ is given by the $d$-separation relations in $\G$~\cite[Section~3.2.2]{Lauritzen_1996}; i.e., a distribution $\p$~is Markov with respect to $\G$ if and only if  $X_A\independent X_B \,\vert\, X_C$ in $\p$ whenever $A$ and $B$ are $d$-separated in $\G$ given $C$.  
The faithfulness assumption asserts that all conditional independence relations entailed by $\p$ are given by $d$-separations in~$\G$~\cite{Spirtes_2001}.  

\begin{assump}[Faithfulness Assumption]
\label{ass: faithfulness}
A distribution $\p$ satisfies the faithfulness assumption with respect to a DAG $\mathcal{G} = ([p],A)$ if for any pair of nodes $i,j\in [p]$ and any $S\subset [p]\backslash\{i,j\}$ we have that
$
i \independent j \, \vert \, S 
$
if and only if $i$ is $d$-separated from $j$ given $S$ in $\mathcal{G}$.
\phantom\qed
\end{assump}

All DAG model learning algorithms assume the Markov assumption, i.e.~the forward direction of the faithfulness assumption, and many of the classical algorithms also assume the converse. 
Unfortunately, the faithfulness assumption has been shown to be restrictive in practice~\cite{URBY13}, and a number of relaxations of this assumption have been suggested~\cite{RZS12}.
For example, restricted faithfulness is the weakest known sufficient condition for consistency of the popular PC-algorithm \cite{Spirtes_2001}.

\begin{assump}[Restricted faithfulness assumption]
\label{ass: restricted faithfulness}
A distribution $\p$ satisfies the restricted faithfulness assumption with respect to a DAG $\mathcal{G} = ([p], A)$ if it satisfies the two conditions:
	\begin{enumerate}
		\item (Adjacency Faithfulness) For all arrows $i\rightarrow j\in A$ we have that
		$X_i\not\independent X_j\,\vert\, X_S$ for all subsets $S\subset[p]\backslash\{i,j\}$.  
		\item (Orientation Faithfulness) For all unshielded triples $(i,j,k)$ and all subsets $S\subset[p]\backslash\{i,k\}$ such that $i$ is $d$-connected to $k$ given $S$, we have that $X_i\not\independent X_k\,\vert\, X_S$.  
		\phantom\qed
	\end{enumerate}
\end{assump}

By sacrificing computation time, some algorithms remain consistent under assumptions that further relax restricted faithfulness, an example being the sparsest permutation algorithm~\cite{RU13}: 
Let $S_p$ denote the space of all permutations of length $p$. 
Given a set of conditional independence relations $\CC$ on $[p]$, every permutation $\pi\in S_p$ is associated to a DAG $\mathcal{G}_\pi$ as follows:
\[
\pi_i\rightarrow\pi_j\in A(\G_\pi) \quad \Leftrightarrow \quad i<j \mbox{ and } \pi_i\not\independent \pi_j \, \vert \, \{\pi_1,\ldots,\pi_{j}\}\backslash\{\pi_i,\pi_j\}.
\]
\begin{figure}[t!]
	\centering
	\begin{tikzpicture}[thick,scale=0.5]
	

	 \node (1) at (0,2) {$1$};
	 \node (2) at (2,0) {$2$};
	 \node (3) at (4,2) {$3$};

	 \node (m1) at (6,2) {$1$};
	 \node (m2) at (8,0) {$2$};
	 \node (m3) at (10,2) {$3$};

	 \node (g1) at (12,2) {$1$};
	 \node (g2) at (14,0) {$2$};
	 \node (g3) at (16,2) {$3$};

	 \draw[->]   (1) -- (2) ;
 	 \draw[->]   (2) -- (3) ;
 	 \draw[->]   (1) -- (3) ;

	 \draw[->]   (m1) -- (m2) ;
 	 \draw[->]   (m3) -- (m2) ;
 	 \draw[->]   (m1) -- (m3) ;

	 \draw[->]   (g1) -- (g2) ;
 	 \draw[->]   (g3) -- (g2) ;

	\node (L1) at (0.5,0) {$\G_\pi$} ; 
	\node (L2) at (6.5,0) {$\G$} ; 
	\node (L3) at (12.5,0) {$\G_\tau$} ; 
	
\draw[fill = blue!25] (18,0) -- (22,0) -- (23,2) -- (22,4) -- (20,4) -- cycle;
 	 \node [circle, draw, fill=black!100, inner sep=2pt, minimum width=2pt] (p1) at (18,0) {};
 	 \node [circle, draw, fill=black!100, inner sep=2pt, minimum width=2pt] (p2) at (22,0) {};
 	 \node [circle, draw, fill=black!100, inner sep=2pt, minimum width=2pt] (p3) at (23,2) {};
 	 \node [circle, draw, fill=black!100, inner sep=2pt, minimum width=2pt] (p4) at (22,4) {};
 	 \node [circle, draw, fill=black!100, inner sep=2pt, minimum width=2pt] (p5) at (20,4) {};

 	 \draw   	 (p1) -- (p2) ;
	 \draw   	 (p2) -- (p3) ;
	 \draw   	 (p3) -- (p4) ;
	 \draw   	 (p4) -- (p5) ;
	 \draw   	 (p5) -- (p1) ;
	 
	 \node (A1) at (16.85,-0.25) {
	 	\begin{tikzpicture}[thick,scale=0.2]
	
			 \node[inner sep =.5pt] (a1) at (10,2) {\scriptsize $1$};
	 		 \node[inner sep =.5pt] (a2) at (12,0) {\scriptsize $2$};
			 \node[inner sep =.5pt] (a3) at (14,2) {\scriptsize $3$};

			 \draw[->]   (a1) -- (a2) ;
 			 \draw[->]   (a3) -- (a2) ;

		\end{tikzpicture}
					   } ;
					   
	 \node (B1) at (19.25,4.5) {
	 	\begin{tikzpicture}[thick,scale=0.2]
	
			 \node[inner sep =.5pt] (b1) at (10,2) {\scriptsize $1$};
	 		 \node[inner sep =.5pt] (b2) at (12,0) {\scriptsize $2$};
			 \node[inner sep =.5pt] (b3) at (14,2) {\scriptsize $3$};

			 \draw[->]   (b1) -- (b2) ;
 			 \draw[->]   (b2) -- (b3) ;
			 \draw[->]   (b1) -- (b3) ;

		\end{tikzpicture}
					   } ;
	\node (B2) at (22.5,4.5) {
	 	\begin{tikzpicture}[thick,scale=0.2]
	
			 \node[inner sep =.5pt] (b1) at (10,2) {\scriptsize $1$};
	 		 \node[inner sep =.5pt] (b2) at (12,0) {\scriptsize $2$};
			 \node[inner sep =.5pt] (b3) at (14,2) {\scriptsize $3$};

			 \draw[->]   (b2) -- (b1) ;
 			 \draw[->]   (b2) -- (b3) ;
			 \draw[->]   (b1) -- (b3) ;

		\end{tikzpicture}
					   } ;
 	 \node (B3) at (23.75,2.25) {
	 	\begin{tikzpicture}[thick,scale=0.2]
	
			 \node[inner sep =.5pt] (b1) at (10,2) {\scriptsize $1$};
	 		 \node[inner sep =.5pt] (b2) at (12,0) {\scriptsize $2$};
			 \node[inner sep =.5pt] (b3) at (14,2) {\scriptsize $3$};

			 \draw[->]   (b2) -- (b1) ;
 			 \draw[->]   (b2) -- (b3) ;
			 \draw[->]   (b3) -- (b1) ;

		\end{tikzpicture}
					   } ;
	
	\node (B4) at (23.2,-0.25) {
	 	\begin{tikzpicture}[thick,scale=0.2]
	
			 \node[inner sep =.5pt] (b1) at (10,2) {\scriptsize $1$};
	 		 \node[inner sep =.5pt] (b2) at (12,0) {\scriptsize $2$};
			 \node[inner sep =.5pt] (b3) at (14,2) {\scriptsize $3$};

			 \draw[->]   (b2) -- (b1) ;
 			 \draw[->]   (b3) -- (b2) ;
			 \draw[->]   (b3) -- (b1) ;

		\end{tikzpicture}
					   } ;

	\end{tikzpicture}
	\caption{For $\CC = \{ 1\independent 3\}$, we see the polytope $\mathcal{A}_3(\CC)$, the graphs $\G_\pi$ and $\G_\tau$ for the $\pi = 123$ and $\tau = 132$, and a graph $\G$ Markov equivalent to $\G_\pi$ that is not a minimal independence map of $\CC$. $\pi$ and $\tau$ are related by transposing $2$ and $3$ in $\pi$ and the arrow $2\rightarrow 3$ in $\G_\pi$ is covered.}
	\label{fig: Minimal I-MAP Examples}
\end{figure}
Examples of the DAGs $\G_\pi$ appear in Figure~\ref{fig: Minimal I-MAP Examples}.
A DAG $\mathcal{G}_\pi$ is known as a minimal independence map with respect to $\CC$, since  any DAG $\mathcal{G}_\pi$ satisfies the minimality assumption with respect to $\CC$, i.e., any conditional independence relation encoded by a $d$-separation in $\mathcal{G}_\pi$ is in $\CC$ and any proper sub-DAG of $\mathcal{G}_\pi$ encodes a conditional independence relation that is not in $\CC$~\cite{Pearl_1988}.  
The sparsest permutation algorithm searches over all DAGs $\mathcal{G}_\pi$ for $\pi\in S_p$ and returns a DAG that maximizes the score
$$
\score(\CC; \mathcal{G}) := 
\begin{cases}
-|\mathcal{G}|	&	\mbox{ if $\mathcal{G}$ is Markov with respect to $\CC$}, \\
-\infty			&	\mbox{otherwise}, \\
\end{cases}
$$
where $|\mathcal{G}|$ denotes the number of arrows in $\mathcal{G}$. 
In \cite{RU13}, it is shown that the sparsest permutation algorithm is consistent under the sparsest Markov representation assumption, which is strictly weaker than restricted faithfulness. 
\begin{assump}[Sparsest Markov representation assumption]
\label{ass: smr}
A probability distribution $\p$ satisfies the sparsest Markov representation assumption with respect to a DAG $\G$ if it is Markov with respect to $\G$ and $|\G|<|\HH|$ for every DAG $\HH$ to which $\p$ is Markov and which satisfies $\HH\notin\M(\G)$.  
\phantom\qed
\end{assump}
The downside to the sparsest permutation algorithm is that it requires a search over all $p!$ permutations of the node set $[p]$.  
A typical approach to accommodate a large search space is to pass to a greedy variant of the algorithm.  
To this end, we analyze greedy variants of the sparsest permutation algorithm in the coming section.

\section{Greedy sparsest permutation algorithm}
\label{sec_greedy_SP}

The sparsest permutation algorithm has a natural interpretation in the setting of discrete geometry.  
The permutohedron on $p$ elements, denoted $\mathcal{A}_p$, is the convex hull in $\R^p$ of all vectors obtained by permuting the coordinates of $(1,2,3,\ldots,p)^T$.  
The sparsest permutation algorithm can be thought of as searching over the vertices of $\mathcal{A}_p$, since it considers the DAGs $\G_\pi$ for each $\pi\in S_p$.  Hence, a natural first step to reduce the size of the search space is to contract all vertices of $\mathcal{A}_p$ that correspond to the same DAG $\mathcal{G}_\pi$. This can be done via the following construction first presented in \cite{MUWY16}.

Two vertices of the permutohedron $\mathcal{A}_p$ are connected by an edge if and only if the permutations indexing the vertices differ by an adjacent transposition.  
We associate a conditional independence relation to adjacent transpositions, and hence to each edge of $\mathcal{A}_p$; namely $\pi_{i}\independent \pi_{i+1}\vert \{\pi_{1}, \dots , \pi_{i-1}\}$ to the edge between
$$
 \quad (\pi_1,\ldots,\pi_i,\pi_{i+1},\ldots,\pi_p)^T \; \mbox{ and } \; (\pi_1,\ldots,\pi_{i+1},\pi_i,\ldots,\pi_p)^T.
$$
In \cite[Section 4]{MUWY16}, it is shown that given a set of conditional independence relations $\mathcal{C}$ from a joint distribution $\mathbb{P}$ on $[p]$, then contracting all edges in $\mathcal{A}_p$ corresponding to conditional independence relations in $\mathcal{C}$ results in a convex polytope, which we denote by  $\mathcal{A}_p(\mathcal{C})$. Note that $\mathcal{A}_p(\emptyset) = \mathcal{A}_p$. Furthermore, if the conditional independence relations in $\mathcal{C}$ form a graphoid, i.e., they satisfy the semigraphoid properties and the intersection property:
\begin{enumerate}
\item if $i\independent j\vert S$ then $j\independent i\vert S$,
\item if $i\independent j\vert S$ and  $i\independent k\vert \{j\}\cup S$, then $i\independent k\vert S$ and $i\independent j\vert \{k\}\cup S$,
\item  if $i\independent j\vert \{k\}\cup S$ and  $i\independent k\vert \{j\}\cup S$, then $i\independent j \vert S$ and $i\independent k\vert S$,
\end{enumerate}
then it was shown in \cite[Theorem~7.1]{MUWY16} that contracting edges in $\mathcal{A}_p$ that correspond to conditional independence relations in $\mathcal{C}$ is the same as identifying vertices of $\mathcal{A}_p$ that correspond to the same DAG. The semigraphoid properties hold for any distribution.  
On the other hand, the intersection property holds, for example, for strictly positive distributions.  
Another example of a graphoid is the set of conditional independence relations $\mathcal{C}$ corresponding to all $d$-separations in a DAG. In that case $\mathcal{A}_p(\mathcal{C})$ is also called a DAG associahedron~\cite{MUWY16}. 
The edge graph of the polytope $\mathcal{A}_p(\mathcal{C})$, where each vertex corresponds to a different DAG, represents a natural search space for a greedy version of the sparsest permutation algorithm.

Through a closer examination of the polytope $\mathcal{A}_p(\mathcal{C})$, we arrive at two greedy versions of the sparsest permutation algorithm: one based on the geometry of $\mathcal{A}_p(\mathcal{C})$ by walking along edges of the polytope and another based on the combinatorial description of the vertices by walking from DAG to DAG. 
These two greedy versions of the sparsest permutation algorithm are given in Algorithms~\ref{alg_edge_SP}~and~\ref{alg_triangle_SP}. 

Both algorithms take as input a set of conditional independence relations $\CC$ and an initial permutation $\pi\in S_p$. 
Beginning at the vertex $\G_\pi$ of $\mathcal{A}_p(\mathcal{C})$, Algorithm~\ref{alg_edge_SP} walks along an edge of $\mathcal{A}_p(\mathcal{C})$ to any vertex whose corresponding DAG has at most as many arrows as $\G_\pi$.  
Once it can no longer discover a sparser DAG, the algorithm returns the last DAG it visited, from which we deduce the corresponding Markov equivalence class.  
Since this algorithm is based on walking along edges of $\mathcal{A}_p(\mathcal{C})$, we call this greedy version the edge sparsest permutation algorithm. 
The corresponding identifiability assumption can be stated as follows:
\begin{assump}[Edge assumption]
\label{ass: edge SP}
A distribution $\p$ satisfies the edge assumption with respect to a DAG $\G$ if it is Markov with respect to $\G$ and if Algorithm~\ref{alg_edge_SP} returns only DAGs in $\M(\G)$.  
\phantom\qed
\end{assump}
{
\begin{algorithm}[t]
\label{alg_edge_SP}
\caption{The Edge Sparsest Permutation Algorithm}
\LinesNumbered
\DontPrintSemicolon
\SetAlgoLined
\SetKwInOut{Input}{Input}
\SetKwInOut{Output}{Output}
\Input{A set of conditional independence relations $\CC$ on node set $[p]$ and a starting permutation $\pi\in S_p$.}
\Output{A minimal independence map $\G$.}
\BlankLine
    Compute the polytope $\mathcal{A}_p(\mathcal{C})$ and set
    $\G := \G_\pi$.
    \;
     Using a depth-first search approach with root $\G$ along the edges of $\mathcal{A}_p(\mathcal{C})$, search for a minimal independence map $\G_\tau$ with $|\G|>|\G_\tau|$.
    If no such $\G_\tau$ exists, return $\G$; else set $\G :=\G_\tau$ and repeat this step. 
\end{algorithm}
}

Algorithm~\ref{alg_edge_SP} requires computing the polytope $\mathcal{A}_p(\mathcal{C})$. This is inefficient, since an edge walk in a polytope only requires knowing the neighbors of a vertex and not the full polytope. In the following, we overcome this inefficiency by providing a graphical characterization of neighboring DAGs. 

We say that an arrow $i\rightarrow j$ in a DAG $\G$ is covered if $\parents(i) = \parents(j)\setminus\{i\}$ and it is  trivially covered if $\parents(i) = \parents(j)\setminus\{i\} =\emptyset$.
For example, the arrows $1\rightarrow 2$ and $2\rightarrow 3$ in the DAG $\G_\pi$ in Figure~\ref{fig: Minimal I-MAP Examples} are both covered, but only the arrow $1\rightarrow 2$ is trivially covered.
In addition, we call a sequence of minimal independence maps $(\G_{\pi^1},\G_{\pi^2},\ldots,\G_{\pi^N})$ a weakly decreasing sequence if $|\G_{\pi^i}|\geq|\G_{\pi^{i+1}}|$ for all $i\in[N-1]$.  
If $\G_{\pi^{i+1}}$ is produced from $\G_{\pi^i}$ by reversing a covered arrow in $\G_{\pi^i}$, then we refer to this sequence as a weakly decreasing sequence determined by covered arrow reversals.  
For instance, given the DAGs $\G_\pi$ and $\G_\tau$ from Figure~\ref{fig: Minimal I-MAP Examples}, $(\G_\pi,\G_\tau)$ is a weakly decreasing sequence determined by covered arrow reversals.
Let $\G_\pi$ and $\G_\tau$ denote two adjacent vertices in a DAG associahedron $\mathcal{A}_p(\mathcal{C})$. Let $\bar{\mathcal{G}}$ denote the skeleton of $\mathcal{G}$; i.e., the undirected graph obtained by undirecting all arrows in $\mathcal{G}$. 
Then, as noted in \cite[Theorem 8.3]{MUWY16}, $\G_\pi$ and $\G_\tau$ differ by a covered arrow reversal if and only if $\overline{\G_\pi}\subseteq \overline{\G_\tau}$ or $\overline{\G_\tau}\subseteq \overline{\G_\pi}$.  
In some instances, this fact gives a combinatorial interpretation of all edges of $\mathcal{A}_p(\CC)$.  
However, this need not always be true as demonstrated in Example~\ref{ex: not all edges correspond to covered arrow flips} in Section~\ref{app: greedy SP} of the Supplementary Material.  

The combinatorial description of some edges of $\mathcal{A}_p(\CC)$ via covered arrow reversals motivates Algorithm~\ref{alg_triangle_SP}, a combinatorial greedy sparsest permutation algorithm. 
Since this algorithm is based on flipping covered arrows, we call this the triangle sparsest permutation algorithm.
Unlike Algorithm~\ref{alg_edge_SP}, this algorithm does not require computing the polytope $\mathcal{A}_p(\CC)$ and is thus the version run in practice. 
Similar to Algorithm~\ref{alg_edge_SP}, we specify an identifiability assumption in relation to Algorithm~\ref{alg_triangle_SP}.  
\begin{assump}[Triangle assumption]
\label{ass: triangle SP}
A distribution $\p$ satisfies the triangle assumption with respect to a DAG $\G$ if it is Markov with respect to $\G$ and if Algorithm~\ref{alg_triangle_SP} returns only DAGs in $\M(\G)$.
\phantom\qed
\end{assump}

{
\begin{algorithm}[t]
\label{alg_triangle_SP}
\caption{The Triangle Sparsest Permutation Algorithm}
\LinesNumbered
\DontPrintSemicolon
\SetAlgoLined
\SetKwInOut{Input}{Input}
\SetKwInOut{Output}{Output}
\Input{A set of conditional independence relations $\CC$ on node set $[p]$ and a starting permutation $\pi\in S_p$.}
\Output{A minimal independence map $\G$.}
\BlankLine
    Set $\G := \G_\pi$.
    \;
    Using a depth-first search approach with root $\G$, search for a minimal independence map $\G_\tau$ with $|\G|>|\G_\tau|$ that is connected to $\G$ by a weakly decreasing sequence determined by covered arrow reversals.
    If no such $\G_\tau$ exists, return $\G$; else set $\G :=\G_\tau$ and repeat this step.
\end{algorithm}
}

In the same way that the sparsest Markov representation assumption is precisely the necessary and sufficient condition under which the sparsest permutation algorithm is consistent, Assumption~\ref{ass: edge SP} and Assumption~\ref{ass: triangle SP}, respectively, are defined to be the necessary and sufficient conditions under which Algorithm~\ref{alg_edge_SP} and Algorithm~\ref{alg_triangle_SP}, respectively, are consistent.  
By associating an identifiability assumption with an algorithm in this way, we can more easily describe which algorithms are consistent for more distributions. 
It is straightforward to verify that every covered arrow reversal in some minimal independence map $\G_\pi$ with respect to $\CC$ corresponds to some edge of the DAG associahedron $\mathcal{A}_p(\CC)$.  
Consequently, if a distribution satisfies the triangle assumption then it also satisfies the edge assumption.
In Theorem~\ref{thm: faithfulness, SMR, and greedy SP assumptions} we show that both these assumptions are weaker than the faithfulness assumption, but stronger than the sparsest Markov representation assumption.

\section{Consistency Guarantees and Identifiability Implications}
\label{sec: greedy sp}

\subsection{Consistency of the edge and triangle sparsest permutation algorithms under faithfulness}
\label{subsec: consistency of greedy sp under faithfulness}
In this section, we prove that both Algorithm~\ref{alg_edge_SP} and Algorithm~\ref{alg_triangle_SP} are pointwise consistent under the faithfulness assumption; i.e., in the oracle-version as $n\to\infty$ the algorithm outputs the true Markov equivalence class.
First note that since the triangle assumption implies the edge assumption, it is sufficient to prove pointwise consistency of Algorithm~\ref{alg_triangle_SP}. 
To prove this, we need to show that for given a set of conditional independence relations $\CC$ corresponding to $d$-separations in a DAG $\mathcal{G}^*$, every weakly decreasing sequence determined by covered arrow reversals ultimately leads to a DAG in $\mathcal{M}(\mathcal{G}^*)$. 
Given two DAGs $\G$ and $\HH$, $\HH$ is an independence map of $\G$, denoted by $\G\leq \HH$, if every conditional independence relation encoded by $\HH$ holds in $\G$ (i.e. $\CI(\G)\supseteq\CI(\HH)$). 
The following simple result, whose proof is given in the Supplementary Material, reveals the main idea of the proof.

\begin{lemma}
\label{thm: faithfulness}
A probability distribution $\p$ on the node set $[p]$ is faithful with respect to a DAG $\G$ if and only if $\G\leq \G_\pi$ for all $\pi\in S_p$.  
\end{lemma}

The goal is to prove that for any pair of DAGs such that $\G_\pi\leq\G_\tau$, there is a weakly decreasing sequence determined by covered arrow reversals such that
$$
\left(\G_\tau=\G_{\pi^0}, \G_{\pi^1}, \G_{\pi^2}, \ldots, \G_{\pi^M}=\G_\pi\right).
$$
Our proof relies heavily on Chickering's consistency proof of greedy equivalence search and, in particular, on his proof of a conjecture known as Meek's conjecture.

\begin{theorem}
\label{thm: meek's conjecture}
\cite[Theorem 4]{C02}
Let $\G$ and $\HH$ be any pair of DAGs such that $\G\leq \HH$.  
Let $r$ be the number of arrows in $\HH$ that have opposite orientation in $\G$, and let $m$ be the number of arrows in $\HH$ that do not exist in either orientation in $\G$.  
There exists a sequence of at most $r+2m$ arrow reversals and additions in $\G$ with the following properties:
\begin{enumerate}
	\item Each arrow reversal is a covered arrow.  
	\item After each reversal and addition, the graph $\G^\prime$ is a DAG and $\G^\prime\leq \HH$.  
	\item After all reversals and additions $\G=\HH$.  
\end{enumerate}
\end{theorem}

In \cite{C02}, a constructive proof of this result is given via the APPLY-EDGE-\\OPERATION algorithm. 
For convenience, we will henceforth refer to this algorithm as the Chickering algorithm.
The Chickering algorithm takes in an independence map $\G\leq \HH$ and adds an arrow to $\G$ or reverses a covered arrow in $\G$ to produce a new DAG $\G^1$ for which $\G\leq \G^1 \leq \HH$.  
By Theorem~\ref{thm: meek's conjecture}, repeated applications of this algorithm produces a sequence of graphs 
$$
\G = \G^0 \leq \G^1 \leq \G^2\leq \cdots \leq \G^N = \HH.  
$$
We will call any sequence of DAGs produced in this fashion a Chickering sequence from $\G$ to $\HH$.  
A quick examination of the Chickering algorithm reveals that there can be multiple Chickering sequences from $\G$ to $\HH$.  
We are interested in identifying a specific type of Chickering sequence in which the covered arrow reversals and edge additions correspond to steps between minimal independence maps in a weakly decreasing sequence. 

Given two DAGs $\G_\pi\leq \G_\tau$, Algorithm~\ref{alg_triangle_SP} proposes that there is a path along the edges of $\mathcal{A}_p(\CC)$ corresponding to covered arrow reversals taking us from $\G_\tau$ to $\G_\pi$, say 
$
\left(\G_{\tau}=\G_{\pi^0}, \G_{\pi^1}, \G_{\pi^2}, \ldots, \G_{\pi^M}=\G_\pi\right),
$
for which $|\G_{\pi^{j-1}}|\geq |\G_{\pi^j}|$ for all $j=1,\ldots,M$.  
Recall that we call such a sequence of minimal independence maps satisfying the latter property a weakly decreasing sequence determined by covered arrow reversals.  
If such a weakly decreasing sequence exists from any $\G_\tau$ to $\G_\pi$, then Algorithm~\ref{alg_triangle_SP} must find it.  
By definition, such a path is composed of covered arrow reversals and arrow deletions.   
Since these are precisely the types of moves used in the Chickering algorithm, then we must understand the subtleties of the relationship between independence maps relating the DAGs $\G_\pi$ for a collection of conditional independence relations $\CC$ and the skeletal structure of the $\G_\pi$.  
To this end, we will use the following two definitions:
We will denote that two DAGs $\G$ and $\HH$ are Markov equivalent by $\G\approx \HH$.
A minimal independence map $\G_\pi$ with respect to a graphoid $\CC$ is called Markov equivalence class-minimal if for all $\G\approx\G_\pi$ and linear extensions $\tau$ of $\G$ we have that $\G_\pi\leq\G_\tau$.  
Notice by \cite[Theorem 8.1]{MUWY16}, it suffices to check only one linear extension $\tau$ for each $\G$.  
The minimal independence map $\G_\pi$ is further called Markov equivalence class-s-minimal if it is class-minimal and $\overline{\G}_\pi\subseteq\overline{\G}_\tau$ for all $\G\approx\G_\pi$ and linear extensions $\tau$ of $\G$.  
We are now ready to state the main proposition that allows us to verify consistency of Algorithm~\ref{alg_triangle_SP} under the faithfulness assumption.  
\begin{proposition}
\label{thm: consistent under faithfulness}
Suppose that $\CC$ is a graphoid and $\G_\pi$ and $\G_\tau$ are minimal independence maps with respect to $\CC$.  
Then
\begin{enumerate}[(i)]
	\item[(a)] if $\G_\pi \approx\G_\tau$ and $\G_\pi$ is class-s-minimal then there exists a weakly decreasing edgewalk from $\G_\pi$ to $\G_\tau$ along $\mathcal{A}_p(\CC)$.  In particular, any Chickering sequence connecting $\G_\tau$ and $\G_\pi$ is a sequence of Markov equivalent minimal independence maps;
	\item[(b)] if $\G_\pi\leq\G_\tau$ but $\G_\pi\not\approx\G_\tau$ then there exists a minimal independence map $\G_{\tau^\prime}$ with respect to $\CC$ satisfying $\G_{\tau^\prime}\leq\G_\tau$ that is strictly sparser than $\G_\tau$ and is connected to $\G_\tau$ by a weakly decreasing edgewalk along $\mathcal{A}_p(\CC)$.  
\end{enumerate}
\end{proposition}

The proof of Proposition~\ref{thm: consistent under faithfulness} can be found in the Supplementary Material.  
We see from Proposition~\ref{thm: consistent under faithfulness}~(a) that class-s-minimality is simply a formality required to guarantee that moving between Markov equivalent minimal independence maps via covered arrow reversals is equivalent to moving along the edge graph of the associahedron $\mathcal{A}_p(\CC)$.  
Recall that, unlike edgewalks along $\mathcal{A}_p(\CC)$, not all Chickering sequences are sequences of minimal independence maps.  
Instead, a single move along an edge of $\mathcal{A}_p(\CC)$ is given by reversing a covered arrow to produce a Markov equivalent graph, taking a linear extension of that graph, and then computing the associated minimal independence map.  
For instance, a single edge of the associahedron depicted in Figure~\ref{fig: Minimal I-MAP Examples} corresponds to transforming the DAG $\G_\pi$ into $\G$ and then into $\G_\tau$.  
Intuitively, by Lemma~\ref{thm: faithfulness} and Theorem~\ref{thm: meek's conjecture}, one would expect that we can always move in such a fashion from a minimal independence map to a sparser one that better approximates the sparsest minimal independence map $\G_{\pi^\ast}$.  
Proposition~\ref{thm: consistent under faithfulness}~(b) says this intuition is correct, and Proposition~\ref{thm: consistent under faithfulness}~(a) ensures that once we make it to the Markov equivalence class of the sparsest minimal independence map, moving between elements of the class is equivalent to moving along the edge graph of $\mathcal{A}_p(\CC)$.  
These ideas form the basis for the proof of the following theorem.
\begin{theorem}
\label{cor: consistent under faithfulness}
Algorithms~\ref{alg_edge_SP} and~\ref{alg_triangle_SP} are pointwise consistent under the faithfulness assumption.
\end{theorem}

\begin{proof}
Since the triangle assumption implies the edge assumption, it suffices to prove consistency of the triangle sparsest permutation algorithm. 
Suppose that $\CC$ is a graphoid that is faithful to the sparsest minimal independence map $\G_{\pi^\ast}$ with respect to $\CC$.  
By Lemma~\ref{thm: faithfulness}, we know that $\G_{\pi^\ast}\leq\G_\pi$ for all $\pi\in S_p$.  
By (b) of Proposition~\ref{thm: consistent under faithfulness}, if Algorithm~\ref{alg_triangle_SP} is at a minimal independence map $\G_\tau$ that is not in the same Markov equivalence class as $\G_\pi^\ast$, then we can take a weakly decreasing edgewalk along $\mathcal{A}_p(\CC)$ to reach a sparser minimal independence map $\G_{\tau^\prime}$ satisfying $\G_{\pi^\ast}\leq\G_{\tau^\prime}\leq\G_\tau$.  
Following repeated applications of Proposition~\ref{thm: consistent under faithfulness}~(b), the algorithm eventually returns a minimal independence map in the Markov equivalence class of $\G_{\pi^\ast}$.
In order for the algorithm to verify that it is in the correct Markov equivalence class, it needs to compute a minimal independence map $\G_\tau$ for a linear extension $\tau$ for each DAG $\G\approx \G_{\pi^\ast}$ and check that it is not sparser than $\G_{\pi^\ast}$.  
Proposition~\ref{thm: consistent under faithfulness}~(a) states that any such minimal independence map $\G_\tau$ is connected to $\G_{\pi^\ast}$ by a weakly decreasing edgewalk along $\mathcal{A}_p(\CC)$.  
In other words, the Markov equivalence class of $\G_{\pi^\ast}$ is a connected subgraph of the edge graph of $\mathcal{A}_p(\CC)$, and hence the algorithm can verify that it has reached the sparsest Markov equivalence class and terminate.  
Thus, Algorithm~\ref{alg_triangle_SP} is pointwise consistent under the faithfulness assumption. 
\end{proof}

\subsection{Consistency of Algorithm~\ref{alg_triangle_SP} using the Bayesian information criterion}
\label{subsec: consistency with BIC}

{
\begin{algorithm}[t]
\label{alg_triangle_SP_BIC}
\caption{Triangle Sparsest Permutation Algorithm with Bayesian information criterion}
\LinesNumbered
\DontPrintSemicolon
\SetAlgoLined
\SetKwInOut{Input}{Input}
\SetKwInOut{Output}{Output}
\Input{Observations $\hat{X}$, initial permutation $\pi$.}
\Output{Permutation $\hat{\pi}$ with DAG $\G_{\hat{\pi}}$.}
\BlankLine
    Set $\hat{\G}_\pi := \underset{\G\, \text{consistent with permutation}\, \pi}{\argmax} \textrm{BIC}(\G; \hat{X})$.
    \;
    Using a depth-first search approach with root $\pi$, search for a permutation $\tau$ with $\textrm{BIC}(\hat{\G}_\tau; \hat{X}) > \textrm{BIC}(\hat{\G}_\pi; \hat{X})$ that is connected to $\pi$ through a sequence of permutations $(\pi_1, \cdots, \pi_k)$ where each permutation $\pi_i$ is produced from $\pi_{i-1}$ by first doing a covered arrow reversal $\hat{\G}_{\pi_{i-1}}$ and selecting a linear extension $\pi_i$ of the DAG $\hat{\G}_{\pi_{i-1}}$.  
    If no such $\hat{\G}_\tau$ exists, return $\hat{\G}_\pi$; else set $\pi :=\tau$ and repeat.
\end{algorithm}
}
We now show that a version of Algorithm~\ref{alg_triangle_SP} that uses the Bayesian information criterion instead of graph sparsity is also consistent under faithfulness.   
This algorithm is  Algorithm~\ref{alg_triangle_SP_BIC}, and it is constructed in analogy to the methods studied in \cite{TK12}.  
\begin{theorem}
\label{thm: alg_bic}
Algorithm~\ref{alg_triangle_SP_BIC} is pointwise consistent under the faithfulness assumption.
\end{theorem}

Theorem~\ref{thm: alg_bic} is proven in Section~\ref{app: proofs of pointwise consistency guarantees for the greedy sparsest permutation algorithm} of the Supplementary Material. 
It is based on the fact that the Bayesian information criterion is locally consistent.  
This fact follows from the first line of the proof of \cite[Lemma 7]{C02}, which states that Bayesian scoring is locally consistent.
We note that Algorithm~\ref{alg_triangle_SP_BIC} differs from the ordering-based search method proposed in  \cite{TK12} in two main ways: First, Algorithm~\ref{alg_triangle_SP_BIC} selects each new permutation by a covered arrow reversal 	in the associated independence maps.  Second, it uses a depth-first-search approach instead of greedy hill-climbing.
In particular, our search guarantees that any independence map of minimal independence maps $\G_\pi\leq \G_\tau$ are connected by a Chickering sequence.  
The proof of Theorem~\ref{thm: alg_bic} then follows since 
$\vert \G_\tau \vert < \vert \G_\pi \vert \textrm{ if and only if }  \textrm{BIC}(\G_\tau; \hat{X}) > \textrm{BIC}(\G_\pi; \hat{X}),$ 
for any minimal independence maps $\G_\pi$ and $\G_\tau$.  
However, since this fact does not hold for arbitrary DAGs satisfying the Markov assumption with respect to a given distribution, the algorithm of \cite{TK12} lacks known consistency guarantees.

\subsection{Beyond faithfulness}
\label{subsec: beyond faithfulness}
We now examine the relationships between the edge, triangle, sparsest Markov representation, faithfulness, and restricted faithfulness assumptions.  
All proofs for this section can be found in Section~\ref{app: proofs of pointwise consistency guarantees for the greedy sparsest permutation algorithm} of the Supplementary Material.
\begin{theorem}
\label{thm: faithfulness, SMR, and greedy SP assumptions}
Faithfulness implies Assumption~\ref{ass: triangle SP}, which implies Assumption~\ref{ass: edge SP}, which implies the sparsest Markov representation assumption.  Moreover, all implications are strict. 
\end{theorem}

It is clear from the definition that restricted faithfulness is a significantly weaker assumption than faithfulness.  
In \cite[Theorem 2.5]{RU13} it was shown that the sparsest Markov representation assumption is strictly weaker than restricted faithfulness. 
In the following, we compare restricted faithfulness to the triangle assumption and show that the restricted faithfulness assumption is not weaker than the triangle assumption. 
\begin{theorem}
\label{thm: TSP and adjacency faithfulness}
Let $\p$ be a semigraphoid  w.r.t.~a DAG $\G$. If $\p$ satisfies the triangle assumption, then $\p$ satisfies adjacency faithfulness with respect to $\G$. However, there exist distributions $\p$ such that $\p$ satisfies the triangle assumption with respect to a DAG $\G$ and $\p$ does not satisfy orientation faithfulness with respect to $\G$.
\end{theorem}

{
	\begin{algorithm}[t!]
		\caption{The Greedy Sparsest Permutation Algorithm} 
		\label{alg: greedy sp depth and start control}
		\LinesNumbered
		\DontPrintSemicolon
		\SetAlgoLined
		\SetKwInOut{Input}{Input}
		\SetKwInOut{Output}{Output}
		\Input{A set of conditional independence relations $\CC$ on node set $[p]$, and two positive integers $d$ and $r$.}
		\Output{A minimal independence map $\G_\pi$.}
		\BlankLine
		Set $R:=0$, and $Y:=\emptyset$. \;
		\While{$R < r$}{ \label{step: select}
			Select a permutation $\pi\in S_n$ and set $\G:=\G_\pi$. 
			\;
			Using a depth-first search approach with root $\G$, search for a minimal independence map $\G_\tau$ with $|\G|>|\G_\tau|$ that is connected to $\G$ by a weakly decreasing sequence determined by covered arrow reversals that is length at most $d$. \label{step_repeat}
			\;  
			\eIf{no such $\G_\tau$ exists}{
				set $Y:=Y\cup\{\G\}$, $R:=R+1$, and go to step~\ref{step: select}.}
			{set $\G :=\G_\tau$ and go to step~\ref{step_repeat}.
			}
		}
		Return the sparsest DAG $\mathcal{G}_\pi$ in the collection $Y$.
	\end{algorithm}
}

\subsection{The problem of Markov equivalence}
\label{subsec: the problem of markov equivalence}
It is important to note that in contrast to for example the PC-algorithm, Algorithms~\ref{alg_edge_SP} and~\ref{alg_triangle_SP} may need to search over DAGs that belong to the same Markov equivalence class. 
This is due to the fact that two DAGs in the same Markov equivalence class differ only by a sequence of covered arrow reversals \cite[Theorem 2]{C95}.  
Thus, the greedy nature of Algorithm~\ref{alg_triangle_SP} can leave us searching through large portions of Markov equivalence classes until we identify a sparser minimal independence map.  
In particular, in order for Algorithm~\ref{alg_triangle_SP} to terminate, it must visit all members of the Markov equivalence class $\M(\G^\ast)$.  

To address this problem, Algorithm~\ref{alg: greedy sp depth and start control} provides a parametrized alternative that approximates Algorithm~\ref{alg_triangle_SP}. 
We call this algorithm the greedy sparsest permutation algorithm since this is the version that we run in practice; see Section~\ref{sec: simulations}.   
The greedy sparsest permutation algorithm operates exactly like Algorithm~\ref{alg_triangle_SP}, with the exception that it bounds the search depth $d$ and number of runs $r$ allowed before the algorithm terminates.  
Recall that Algorithm~\ref{alg_triangle_SP} searches for a weakly decreasing edge-walk from a minimal independence map $\G_\pi$ to another $\G_\tau$ with $|\G_\pi|>|\G_\tau|$ via a depth-first-search approach.  
In Algorithm~\ref{alg: greedy sp depth and start control}, if this search step does not produce a sparser minimal independence map after searching up to and including depth $d$, the algorithm terminates and returns $\G_\pi$.  
The computational analysis in~\cite{GP01} suggests that the average Markov equivalence class contains four graphs.  
This suggests that a search depth of $4$ is, on average, sufficient for escaping a Markov equivalence class of minimal independence maps.  
This intuition is verified via simulations in Section~\ref{sec: simulations}.

\section{Uniform Consistency}
\label{sec: uniform consistency}

 In this section we show that minor adjustments can turn Algorithm~\ref{alg_triangle_SP} into Algorithm~\ref{alg:grsp}, which is uniformly consistent in the high-dimensional Gaussian setting.  
 In particular, Algorithm~\ref{alg:grsp} only tests conditioning sets made up of parent nodes of covered arrows.
 This feature turns out to be critical for high-dimensional consistency.  
Recently, it was shown that a variant of greedy equivalence search is consistent in the high-dimensional setting~\cite{Marloes16}. 
Similarly to this approach, by assuming sparsity of the initial DAG, we obtain uniform consistency of the triangle sparsest permutation algorithm in the high-dimensional setting, i.e., it converges to the data-generating DAG when $p$ scales with $n$. 

\begin{algorithm}[t]
\caption{A High-dimensional Greedy Sparsest Permutation Algorithm}
\label{alg:grsp}
\LinesNumbered
\KwInput{Observations $\hat{X}$, threshold $\lambda$, and initial permutation $\pi_0$.}
\KwOutput{Permutation $\hat{\pi}$ together with the DAG $\hat{\mathcal{G}}_{\hat{\pi}}$.}
Construct the minimal independence map $\hat{\mathcal{G}}_{\pi_0}$ from the initial permutation $\pi_0$ and $\hat{X}$\;
Perform Algorithm~\ref{alg_triangle_SP} with constrained conditioning sets, i.e., let $i \rightarrow j$ be a covered arrow and let $S=\textrm{pa}(i)=\textrm{pa}(j)\setminus\{i\}$; perform the edge flip, i.e.~$i \leftarrow j$, and update the DAG by removing edges $(k,i)$ for $k\in S$ such that $\vert \hat{\rho}_{i, k \vert (S \union \lbrace j \rbrace \setminus \lbrace k \rbrace)} \vert \leq \lambda$ and edges $(k,j)$ for $k\in S$ such that  $\vert \hat{\rho}_{j, k \vert (S\setminus \lbrace k \rbrace)} \vert \leq \lambda$.
\end{algorithm}

Letting the dimension $p$ grow as a function of the sample size $n$, we write $p=p_n$.  
Similarly, for the true underlying DAG and the data-generating distribution we let $G^*=G^*_n$ and $\PP=\PP_n$, respectively. 
The assumptions under which we will guarantee high-dimensional consistency of Algorithm~\ref{alg:grsp} are as follows:
\begin{enumerate}
\item[(1)] $\PP_n$ is multivariate Gaussian and faithful to the DAG $\mathcal{G}^*_n$ for all $n$.
\item[(2)] The number of nodes $p_n$ scales as $p_n = \mathcal{O}(n^a)$ for some $0 \leq a < 1$.
\item[(3)] Given an initial permutation $\pi_0$, the maximal degree $d_{\pi_0}$ of the corresponding minimal independence map $\mathcal{G}_{\pi_0}$ satisfies $d_{\pi_0} = \mathcal{O}(n^{1 - m})$ for some $0 < m \leq 1$.
\item[(4)] There exists $M<1$ and $c_n > 0$ such that all non-zero partial correlations $\rho_{i,j \vert S}$ satisfy $\vert \rho_{i,j \vert S} \vert \leq M$ and $\vert \rho_{i,j \vert S} \vert \geq c_n$ where $c_n^{-1} = \mathcal{O}(n^\ell)$ for some $0 < \ell < m / 2$.
\end{enumerate}
Analogous to the conditions needed in \cite{KB07}, assumptions (1), (2), (3), and (4) relate to faithfulness, the scaling of the number of nodes with the number of observations, the maximum degree of the initial DAG, and bounds on the minimal non-zero and maximal partial correlations, respectively.  
In the Gaussian setting, the conditional independence relation $X_j  \independent X_k \vert X_S$ is equivalent to the partial correlation $\rho_{j,k\vert S} = \textrm{corr}(X_j,X_k \vert X_S)$ equaling zero, and a hypothesis test based on Fischer's $z$-transform can be used to test whether $X_j  \independent X_k \vert X_S$. 
Combining these facts, we arrive at the following theorem.

\begin{theorem}
\label{thm_high_dim_consist}
Suppose that assumptions~(1), (2), (3), and (4) hold and let the threshold $\lambda$ in Algorithm~\ref{alg:grsp} be defined as $\lambda := c_n /2$.  
Then there exists a constant $c>0$ such that Algorithm~\ref{alg:grsp} is consistent, i.e., it returns a DAG $\hat{\mathcal{G}}_{\hat\pi}$ that is in the same Markov equivalence class as $\mathcal{G}^*_n$, with probability at least $1 - \mathcal{O}\{\exp(-cn^{1 - 2\ell})\}$, where $\ell$ is defined to satisfy assumption (4).
\end{theorem}

As seen in the proof of Theorem~\ref{thm_high_dim_consist}, consistent estimation in the high-dimensional setting requires that we initialize the algorithm at a permutation satisfying assumption (3). This assumption corresponds to a sparsity constraint. In the Gaussian oracle setting the problem of finding a sparsest DAG is equivalent to finding the sparsest Cholesky decomposition of the inverse covariance matrix~\cite{RU13}.  Various heuristics have been developed for finding sparse Cholesky decompositions, the most prominent being the minimum degree algorithm~\cite{George_1973, TW67}. In Algorithm~\ref{alg:emmd} we provide a heuristic for finding a sparse minimal independence map $\mathcal{G}_\pi$ that reduces to the minimum degree algorithm in the oracle setting as shown in Theorem~\ref{thm: NBMDA consistent in non-oracle setting}.
In Algorithm~\ref{alg:emmd}, for a subset of nodes $S\subset[p]$ we let $G_S$ denote the vertex-induced subgraph of $G$ with node set $S$, and for $k\in V$ we let $\adj(G, i)$ denote the nodes $k\in V\setminus\{i\}$ such that $\{i,k\}\in E$. 
The following theorem states that Algorithm~\ref{alg:emmd} is equivalent to the minimum degree algorithm \cite{TW67} in the oracle setting.  
\begin{algorithm}[t]
	\caption{A neighbor-based minimum degree algorithm}
	\label{alg:emmd}
	\KwInput{Observations $\hat{X}$, threshold $\lambda$}
	\KwOutput{Permutation $\hat{\pi}$ together with the DAG $\hat{\mathcal{G}}_{\hat{\pi}}$}
	Set $S := [p]$; construct undirected graph $\hat{G}_S$ with $(i,j) \in \hat{G}_S$ if and only if $\vert\hat{\rho}_{i,j \vert (S \setminus \lbrace i, j \rbrace)}\vert \geq \lambda$\;
	\While{$S \neq \emptyset$}{
		Uniformly draw node $k$ from all nodes with the lowest degree in the graph $\hat{G}_S$\;
		Construct $\hat{G}_{S \setminus \lbrace k \rbrace}$ by first removing node $k$ and its adjacent edges; then update the graph $\hat{G}_{S \setminus \lbrace k \rbrace}$ as follows:
		\begin{align*}
		\begin{split}
		\textrm{$\forall i, j \in \adj(\hat{G}_S, k)$:}\, & \textrm{if } (i,j) \textrm{ not an edge in } \hat{G}_S, \textrm{add } (i,j); \\
		& \textrm{else } (i,j) \textrm{ an edge in } \hat{G}_{S \setminus \lbrace k \rbrace}\; \textrm{iff}\; \vert\hat{\rho}_{i,j \vert S \setminus \lbrace i, j, k \rbrace}\vert \geq \lambda; \\
		\textrm{$\forall i, j \notin \adj(\hat{G}_S, k)$:}\, & (i,j) \textrm{ an edge in } \hat{G}_{S \setminus \lbrace k \rbrace}\; \textrm{iff}\; (i,j) \textrm{ an edge in } \hat{G}_S. \\
		\end{split}
		\end{align*}
		Set $\hat{\pi}(k) := \vert S \vert$ and $S := S \setminus \lbrace k \rbrace$.
	}
	Output the minimal independence map $\hat{\mathcal{G}}_{\hat{\pi}}$ constructed from $\hat{\pi}$ and $\hat{X}$.
\end{algorithm}

\begin{theorem}
\label{thm: NBMDA is equivalent to MDA in oracle setting}
Let the data-generating distribution $\PP$ be multivariate Gaussian with precision matrix $\Theta$.  
Then in the oracle-setting the set of possible output permutations from Algorithm~\ref{alg:emmd} is equal to the possible output permutations of the minimum degree algorithm applied to~$\Theta$.
\end{theorem}
The following result shows that Algorithm~\ref{alg:emmd} in the non-oracle setting is also equivalent to the minimum degree algorithm in the oracle setting.
\begin{theorem}
\label{thm: NBMDA consistent in non-oracle setting}
Suppose that assumptions (1), (2), and (4) hold, and let the threshold $\lambda$ in Algorithm~\ref{alg:grsp} be defined as $\lambda := c_n /2$.  Then with probability at least $1 - \mathcal{O}\{\exp(-cn^{1-2\ell})\}$ the output permutation from Algorithm~\ref{alg:emmd} is contained in the possible output permutations of the minimum degree algorithm applied to $\Theta$.
\end{theorem}

\section{Simulations}
\label{sec: simulations}

The simulations presented in this section were done using the \verb+R+ library \verb+pcalg+~\cite{pcalg}, and linear structural equation models with Gaussian noise:
$$
(X_1,\ldots,X_p)^T = \{(X_1,\ldots,X_p)A\}^T +\varepsilon,
$$
where $\varepsilon\sim\mathcal{N}(0, \mathbb{I}_p)$ with $\mathbb{I}_p$ being the $p\times p$ identity matrix and $A = [a_{ij}]_{i,j=1}^p$ is, without loss of generality, an upper-triangular matrix of edge weights with $a_{ij}\neq 0$ if and only if $i\rightarrow j$ is an arrow in the underlying DAG $\G^{*}$. 
For each simulation study, we generated 100 realizations of a $p$-node random Gaussian DAG model on an Erd\"os-Renyi graph for different values of $p$ and expected neighborhood sizes; i.e., edge probabilities. 
The edge weights $a_{ij}$ were sampled uniformly in $[-1,-0.25]\cup [0.25, 1]$, ensuring that the edge weights are bounded away from 0. 
We first analyzed the oracle setting, where we have access to the true underlying covariance matrix $\Sigma$. 
In the remaining simulations, $n$ samples were drawn from the distribution induced by the Gaussian DAG model for different values of $n$ and $p$. 
In the oracle setting, the conditional independence relations were computed by thresholding the partial correlations using different thresholds $\lambda$.   
For the simulations with $n$ samples, conditional independence relations were estimated by applying Fisher's $z$-transform and comparing the derived $p$-values with a significance level $\alpha$. 
In the oracle and low-dimensional settings, the greedy equivalence search, denoted GES in all figures, was simulated using the Bayesian information criterion.  
In the high-dimensional setting, we used the $\ell_0$-penalized maximum likelihood estimation score~\cite{Marloes16,VAN13}. 

\begin{figure}[t!]
	\centering
	\subfigure[$p=10$, $\lambda=0.001$]{\includegraphics[width=0.35\textwidth]{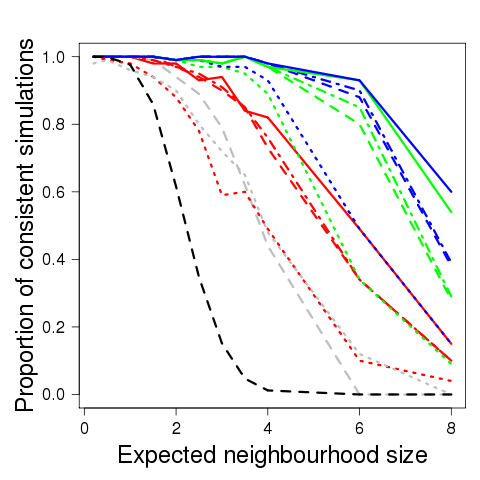}\label{fig:10_0001}}
	\subfigure[$p=10$, $\lambda=0.01$]{\includegraphics[width=0.35\textwidth]{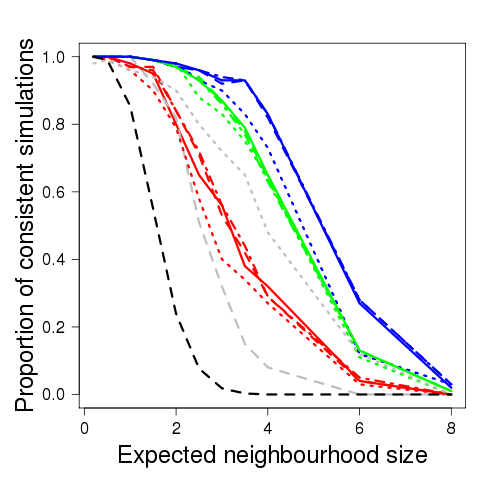}\label{fig:10_001}}
	\subfigure[$p=10$, $\lambda=0.1$]{\includegraphics[width=0.35\textwidth]{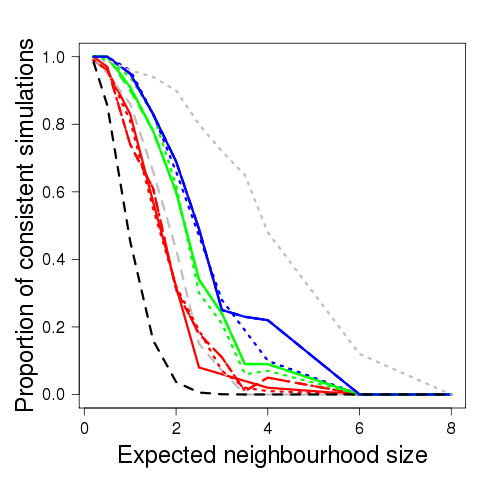}\label{fig:10_01}}\hspace{18pt}
	\subfigure[Legend]{\includegraphics[width=0.3\textwidth]{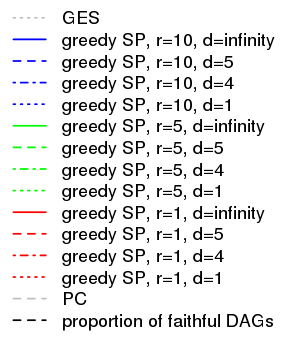}\label{fig:legend}}
	\caption{Expected neighborhood size versus proportion of consistently recovered Markov equivalence classes based on 100 simulations for each expected neighborhood size on DAGs with $p=10$ nodes, edge weights sampled uniformly in $[-1,-0.25]\cup [0.25,1]$, and $\lambda$-values $0.1$, $0.01$ and $0.001$. Greedy SP denotes Algorithm~\ref{alg: greedy sp depth and start control}.  When $r = 1$ and $d = \infty$ this is Algorithm~\ref{alg_triangle_SP}.}
	\label{fig: proportion consistent DAGs}
\end{figure} 

\begin{figure}[t!]
\centering
\subfigure[$s = 0.2$, $n=300$]{\includegraphics[width=0.32\textwidth]{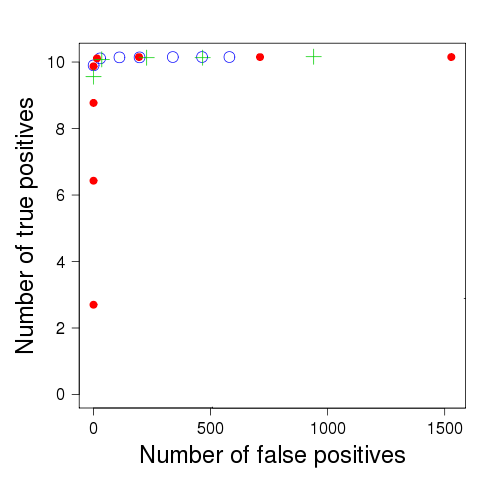}}
\subfigure[$s = 1$, $n=300$]{\includegraphics[width=0.32\textwidth]{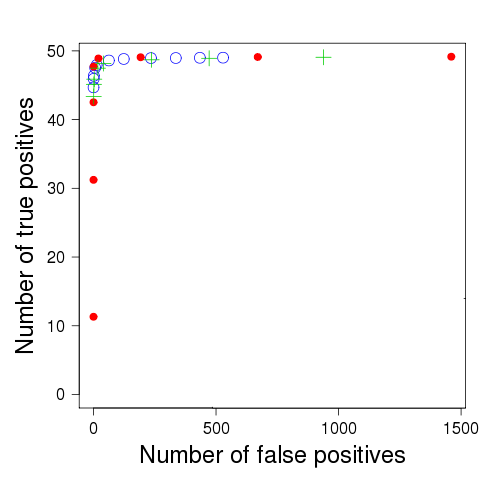}}
\subfigure[$s = 2$, $n=300$]{\includegraphics[width=0.32\textwidth]{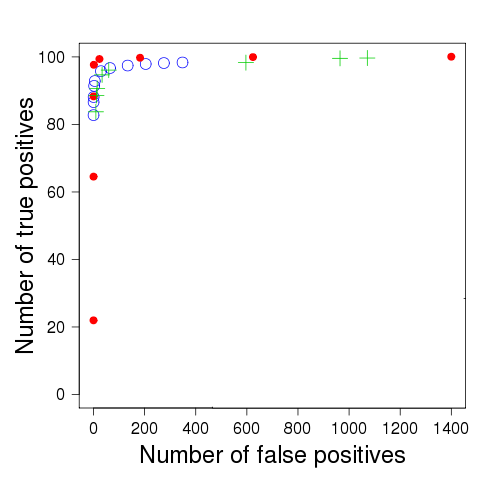}}
\subfigure[$s = 0.2$, $n=300$]{\includegraphics[width=0.32\textwidth]{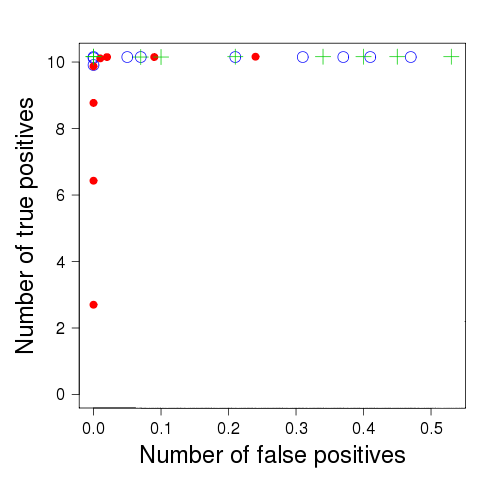}}
\subfigure[$s = 1$, $n=300$]{\includegraphics[width=0.32\textwidth]{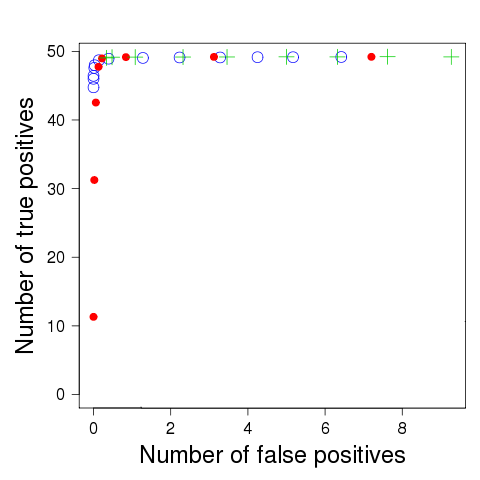}}
\subfigure[$s = 2$, $n=300$]{\includegraphics[width=0.32\textwidth]{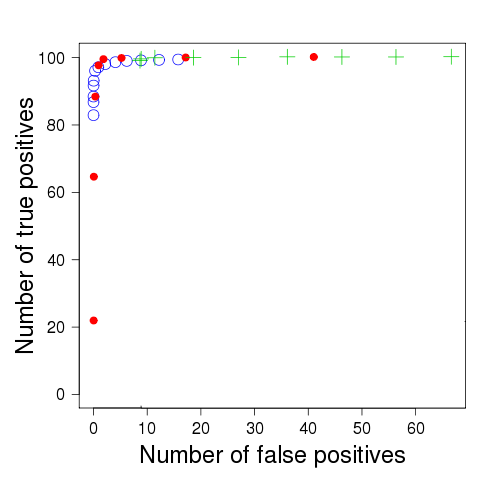}}
\vspace{-0.2cm}
\caption{ROC curves for skeleton recovery for 100 simulations on DAGs with $100$ nodes, expected neighborhood size $s$, sample size $n$, and edge weights sampled uniformly in $[-1,-0.25]\cup [0.25,1]$.  Figures (a)-(c) are without prior knowledge of the moral graph(d)-(f) are with prior knowledge of the moral graph.  Dots denote high-dim GES, crosses denote high-dim greedy SP, and circles denote high-dim PC.}
	\label{fig:high_dim_roc}
	\vspace{-0.4cm}
\end{figure}

Figure~\ref{fig: proportion consistent DAGs} compares the proportion of consistently estimated DAGs in the oracle setting for Algorithm~\ref{alg: greedy sp depth and start control} with number of runs $r \in \lbrace 1, 5, 10\rbrace$ and depth $d\in\{1,4,5,\infty\}$, the PC-algorithm, and the greedy equivalence search. 
Notice that the instance of Algorithm~\ref{alg: greedy sp depth and start control} with parameter settings $r = 1$ and $d = \infty$ is the triangle sparsest permutation algorithm; i.e.~Algorithm~\ref{alg_triangle_SP}.
The number of nodes in these simulations is $p=10$, and we consider $\lambda$-values: $0.1, 0.01$ and $0.001$ for the PC-algorithm and Algorithm~\ref{alg: greedy sp depth and start control}. 
Note that we only run the greedy equivalence search with $n=100,000$ samples since there is no oracle version for this algorithm.
As expected, increasing the number of runs for Algorithm~\ref{alg: greedy sp depth and start control} results in a consistently higher rate of model recovery. In addition, for each fixed number of runs, Algorithm~\ref{alg: greedy sp depth and start control} with search depth $d=4$ performs similarly to $d=\infty$, in line with the observation that the average Markov equivalence class has $4$ elements, as discussed in Section~\ref{subsec: the problem of markov equivalence}.  
For this reason, we recommend setting the search depth $d = 4$. Regarding the choice of $r$, in the low-dimensional setting we have found that choosing $r$ to be of the same magnitude as the number of nodes $p$ produces good estimates. To accelerate computations, for high-dimensional sparse graphs with large $p$, we used $d=1$ and $r=50$; see Figure~\ref{fig:high_dim_roc}.

\begin{figure}[!t]
	\centering
	\subfigure[$p=10$, $\lambda=0.001$]{\includegraphics[width=0.35\textwidth]{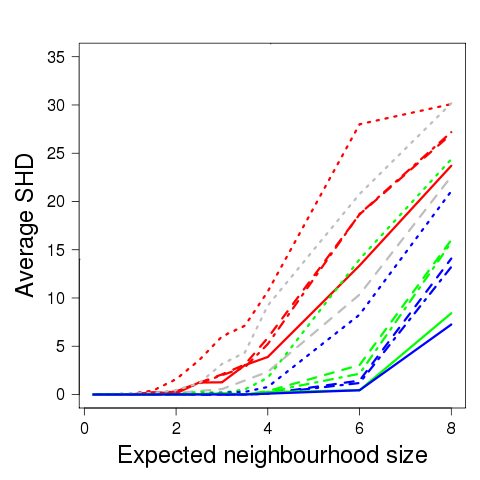}\label{fig:10_0001_shd}}
	\subfigure[$p=10$, $\lambda=0.01$]{\includegraphics[width=0.35\textwidth]{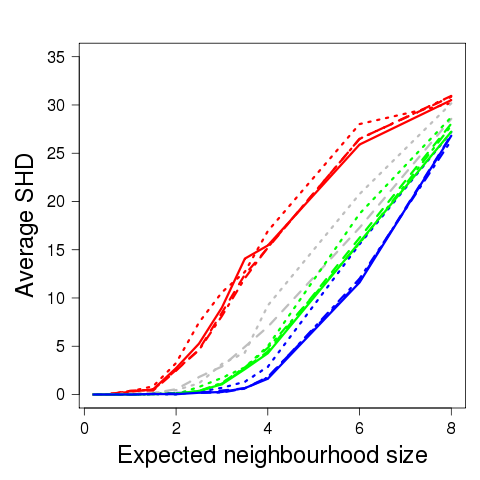}\label{fig:10_001_shd}}
	\subfigure[$p=10$, $\lambda=0.1$]{\includegraphics[width=0.35\textwidth]{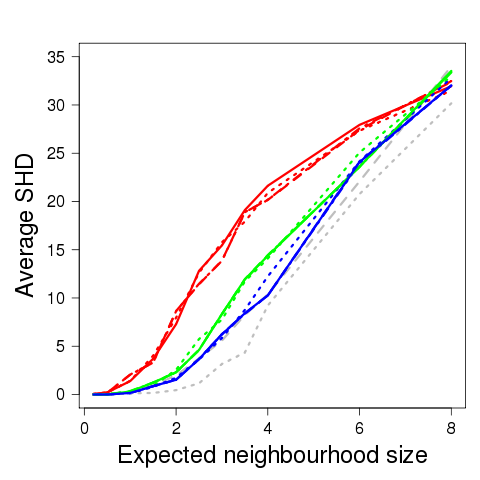}\label{fig:10_01_shd}}\hspace{18pt}
	\subfigure[Legend]{\includegraphics[width=0.3\textwidth]{legend.png}\label{fig:legend}}
	\vspace{-0.2cm}
	\caption{Expected neighborhood size versus structural Hamming distance between the true and recovered Markov equivalence classes based on 100 simulations for each expected neighborhood size on DAGs with $p=10$ nodes, edge weights sampled uniformly in $[-1,-0.25]\cup [0.25,1]$, and $\lambda$-values $0.1$, $0.01$ and $0.001$.}
	\label{fig: structural hamming distance}
	\vspace{-0.2cm}
\end{figure}

For each run, we also recorded the structural Hamming distance between the true and the recovered Markov equivalence classes. 
Figure~\ref{fig: structural hamming distance} shows the average structural Hamming distance versus the expected neighborhood size of the true DAG. 
While Figure~\ref{fig: proportion consistent DAGs} demonstrates that Algorithm~\ref{alg: greedy sp depth and start control}  with search depth $d=4$ and multiple runs learns the true Markov equivalence class at a higher rate than the PC-algorithm and greedy equivalence search when $\lambda$ is chosen small, Figure~\ref{fig: structural hamming distance} shows that, for small values of $d$ and $r$, when Algorithm~\ref{alg: greedy sp depth and start control} learns the wrong DAG it is further off from the true DAG than the PC-algorithm. On the other hand, it appears that this trend only holds for Algorithm~\ref{alg: greedy sp depth and start control}  with a relatively small search depth and few runs. 
That is, increasing the value of these parameters ensures that the wrong DAG learned by Algorithm~\ref{alg: greedy sp depth and start control}  will consistently be closer to the true DAG than that learned by the PC-algorithm.  

\begin{figure}[t!]
	\vspace{-0.5cm}
	\centering
	\subfigure[Directed edge; $n = 1000$]{\includegraphics[width=0.35\textwidth]{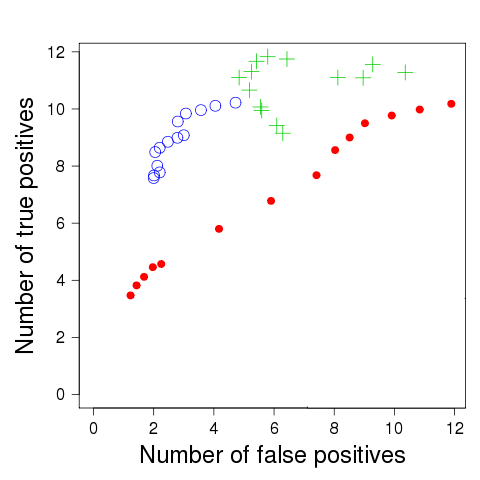}\label{fig:10_0001}}
	\subfigure[Skeleton; $n = 1000$]{\includegraphics[width=0.35\textwidth]{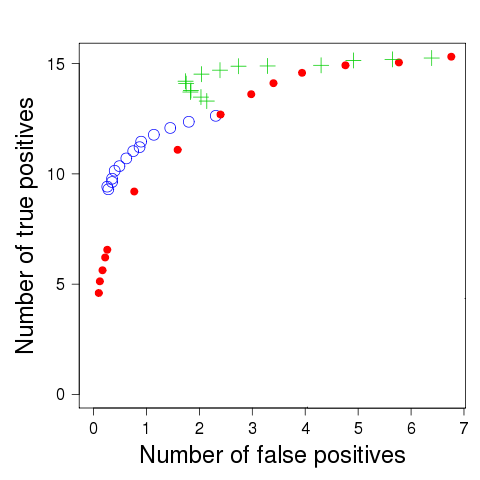}\label{fig:10_001}}
	\subfigure[Directed edge; $n = 10000$]{\includegraphics[width=0.35\textwidth]{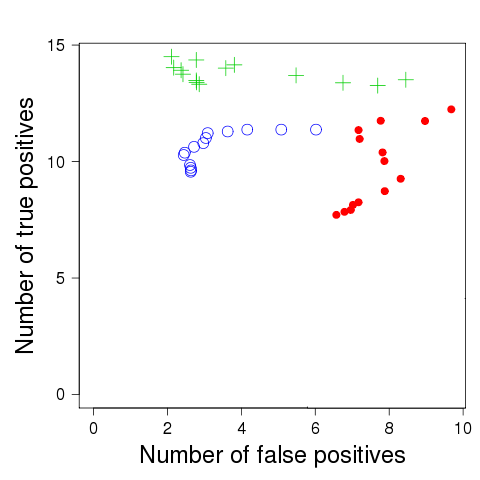}\label{fig:10_01}}
	\subfigure[Skeleton; $n = 10000$]{\includegraphics[width=0.35\textwidth]{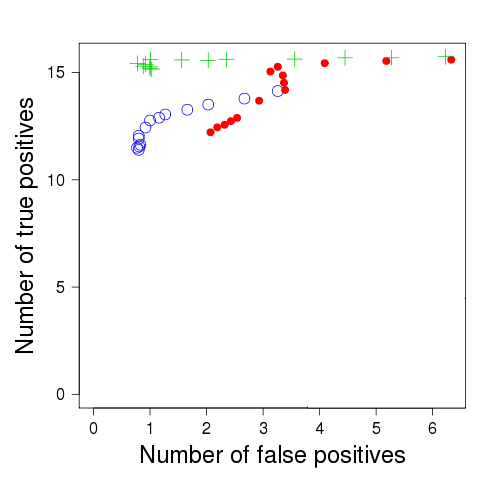}\label{fig:legend}}
	\vspace{-0.1cm}
	\caption{Receiver operating characteristic curves for directed edge recovery and skeleton recovery based on 100 simulations on DAGs with $8$ nodes, expected neighborhood size $4$ and sample size $n\in\{1000,10000\}$. The dots denote GES, crosses is greedy SP, and circles is PC.}
	\label{fig: low dimensional ROCs}
\end{figure}

Recall that Algorithm~\ref{alg: greedy sp depth and start control} and the PC-algorithm can be sensitive to wrong conditional independence test results.  
Each edge and non-edge in the DAG returned by Algorithm~\ref{alg: greedy sp depth and start control} or the PC-algorithm is the result of a conditional independence test.  
To get a sense of their respective sensitivities to wrong conditional independence tests, in Figure~\ref{fig: low dimensional ROCs} we report the number of true positives and false positives for directed edge recovery and skeleton recovery.  
Each data point in the plots represent the average number of true positives and false positives based on $100$ simulated models with $p = 8$ nodes and expected neighborhood size $4$ with a fixed parameter setting.
For Algorithm~\ref{alg: greedy sp depth and start control} and the PC-algorithm, the reported data points correspond to $14$ chosen significance levels in the interval $[0.00005,0.6]$.  
In practice, we recommend tuning the significance level parameter via stability selection~\cite{MB10} as described in Section~\ref{sec_real data}. 
Similarly, for greedy equivalence search the reported data points correspond to $14$ different choices of the scaling constant $c$ from the interval $[0.125, 100]$ for the $\ell_0$-penalization parameter $\lambda_n = c \log (n)$. In Figure~\ref{fig: low dimensional ROCs}, we see that Algorithm~\ref{alg: greedy sp depth and start control} generally outperforms greedy equivalence search and the PC-algorithm in both directed edge and skeleton recovery with large enough sample size.

As noted in Section~\ref{sec_greedy_SP}, we do not consider Algorithm~\ref{alg_edge_SP} in these simulations. 
Recall that Algorithm~\ref{alg_triangle_SP} relies on the fact that a, generally strict, subset of the edges of the polytope $\mathcal{A}_p(\CC)$ have a combinatorial interpretation in terms of their associated minimal independence maps, which makes moving between elements of the search space easier to code.  
Since we do not have a complete combinatorial characterization of all edges of $\mathcal{A}_p(\CC)$, implementing Algorithm~\ref{alg_edge_SP} would require generating a geometric realization of $\mathcal{A}_p(\CC)$ in a program such as \verb+polymake+~\cite{GJ97}, recovering the complete edge graph of this embedding, and then implementing our search over this data structure.  
A natural line of follow-up research is to identify a complete combinatorial interpretation of the edges of $\mathcal{A}_p(\CC)$ so as to allow for an implementation of Algorithm~\ref{alg_edge_SP} that does not require computing the entire polytope $\mathcal{A}_p(\CC)$ and its edge graph.
As shown in Theorem~\ref{thm: faithfulness, SMR, and greedy SP assumptions}, an efficient implementation of Algorithm~\ref{alg_edge_SP} should recover the true DAG at a higher rate than Algorithm~\ref{alg_triangle_SP}.
\begin{figure}[t!]
	\centering
	\subfigure[$n = 1,000$, $\alpha = 0.0001$]{\includegraphics[width=0.3\textwidth]{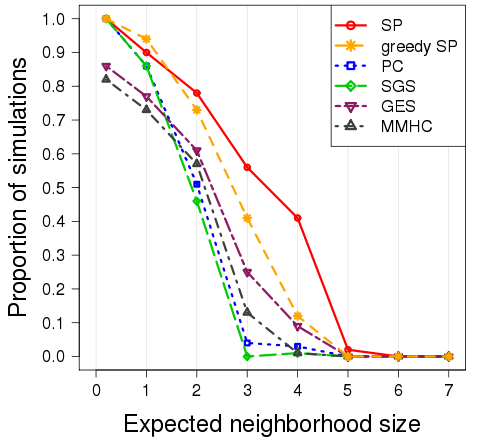}}
	\subfigure[$n = 1,000$, $\alpha = 0.001$]{\includegraphics[width=0.3\textwidth]{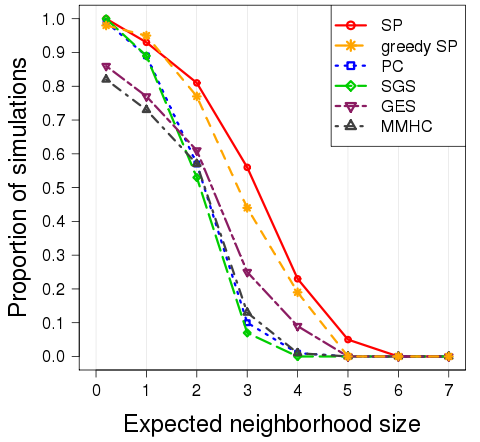}}
		\subfigure[$n = 1,000$, $\alpha = 0.01$]{\includegraphics[width=0.3\textwidth]{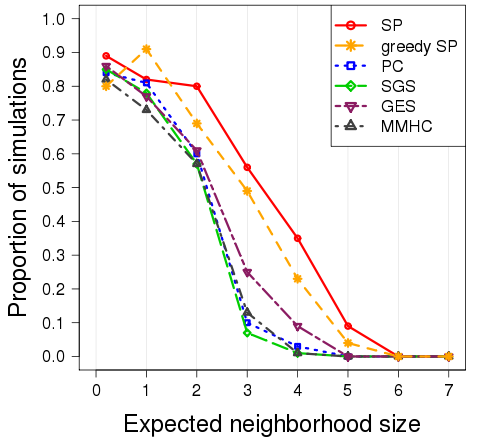}}
\subfigure[$n = 10,000$, $\alpha = 0.0001$]{\includegraphics[width=0.3\textwidth]{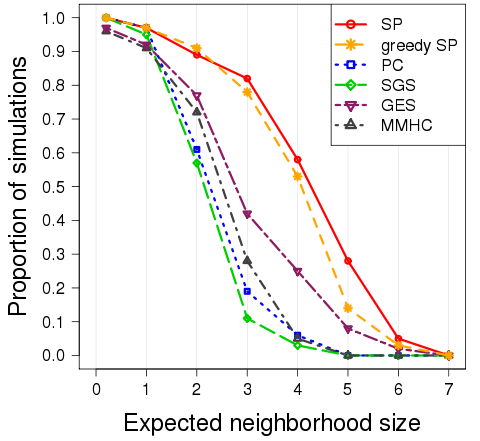}}
	\subfigure[$n = 10,000$, $\alpha = 0.001$]{\includegraphics[width=0.3\textwidth]{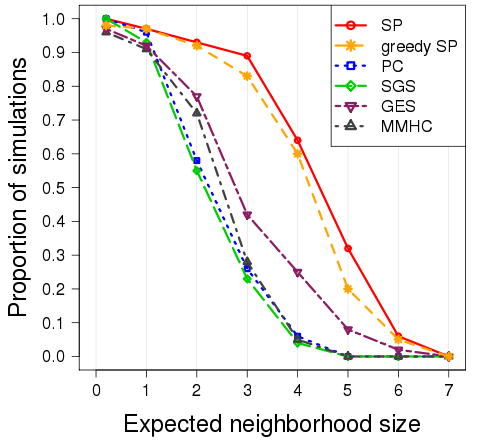}}
		\subfigure[$n = 10,000$, $\alpha = 0.01$]{\includegraphics[width=0.3\textwidth]{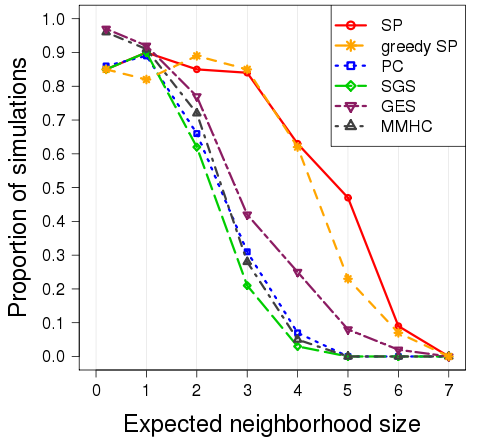}}
		\vspace{-0.2cm}
	\caption{Expected neighborhood size versus proportion of consistently recovered skeleta based on 100 simulations for each expected neighborhood size on DAGs with $p=8$ nodes, sample size $n = 1,000$ and $10,000$, edge weights sampled uniformly in $[-1,-0.25]\cup [0.25,1]$, and $\alpha$-values $0.01$, $0.001$ and $0.0001$;  we used $r =10$ and $d = 4$ for the greedy sparsest permutation algorithm.}
	\vspace{-0.5cm}
	\label{fig:low_dim_consist}
\end{figure}

We then compared the recovery performance of Algorithm~\ref{alg: greedy sp depth and start control} to the sparsest permutation algorithm, greedy equivalence search, the PC-algorithm and its original version, denoted SGS in Figure~\ref{fig:low_dim_consist}, and the max-min hill-climbing algorithm~\cite{TBA2006}, which is denoted MMHC in Figure~\ref{fig:low_dim_consist}. 
This hybrid method first estimates a skeleton through conditional independence testing and then performs a hill-climbing search to orient the edges.
We fixed the number of nodes to be $p=8$ due to the computational limitations of the sparsest permutation algorithm, and considered sample sizes $n= \lbrace 1,000, 10,000\rbrace$. 
We analyzed the performance of greedy equivalence search using the Bayesian information criterion along with Algorithm~\ref{alg: greedy sp depth and start control} and the PC-algorithm for $\alpha = \lbrace 0.01, 0.001, 0.0001 \rbrace$. 
Figure~\ref{fig:low_dim_consist} shows that the sparsest permutation and greedy sparsest permutation algorithms achieve the best performance among all algorithms. 
Since for computational reasons the sparsest permutation algorithm cannot be applied to graphs with over 10 nodes, Algorithm~\ref{alg: greedy sp depth and start control} is the preferable approach for most applications.

In the remainder of this section, we analyze the performance of Algorithm~\ref{alg:grsp} in the sparse high-dimensional setting.  
We compared the performance of Algorithm~\ref{alg:grsp} with $d = 1$ and $r = 50$ with methods that have high-dimensional consistency guarantees; namely the PC-algorithm~\cite{KB07} and greedy equivalence search \cite{Marloes16,VAN13}.
The initial permutation of Algorithm~\ref{alg:emmd} and its associated minimal independence map were used as a starting point in Algorithm~\ref{alg:grsp}, called high-dim greedy SP in Figure~\ref{fig:high_dim_roc}.  
To understand the influence of accurately selecting an initial minimal independence map on the performance of Algorithm~\ref{alg:grsp}, we also considered the case when the moral graph of the data-generating DAG is given as prior knowledge; these results appear in Figure~\ref{fig:high_dim_roc} (d)-(f).  
Figure~\ref{fig:high_dim_roc} compares the skeleton recovery of Algorithm~\ref{alg:grsp} with the PC-algorithm and greedy equivalence search, both without prior knowledge of the moral graph (subfigures (a)-(c)), and with prior knowledge of the moral graph (subfigures (d)-(f)). 
We used the ARGES-CIG algorithm~\cite{Marloes16} to run greedy equivalence search with knowledge of the moral graph.

The number of nodes in our simulations is $p=100$, the number of samples considered is $n = 300$, and the neighborhood sizes used are $s = 0.2$, $1$ and $2$.
We varied the tuning parameters of each algorithm; namely, the significance level $\alpha$ for the PC-algorithm and Algorithm~\ref{alg:grsp}, and the penalization parameter $\lambda_n$ for greedy equivalence search.  
We reported the average number of true positives and false positives for each tuning parameter in the plots shown in Figure~\ref{fig:high_dim_roc}.
This figure shows that, unlike the low-dimensional setting, although Algorithm~\ref{alg:grsp} is still comparable to the PC-algorithm and greedy equivalence search in the high-dimensional setting, greedy equivalence search tends to achieve a slightly better performance in some of the settings.

\section{An Application to Real Data}
\label{sec_real data}
In this section, we compare the performance of the greedy sparsest permutation algorithm, i.e., Algorithm~\ref{alg: greedy sp depth and start control}, with that of the PC-algorithm and greedy equivalence search on the task of gene regulatory network recovery. We consider the perturb-seq data set \cite{D16} containing both observational and interventional data from bone-marrow derived dendritic cells.  
Each data point contains gene expression measurements of 32,777 genes; each interventional data point is sampled from a cell where a single gene was targeted for deletion using the CRISPR/Cas9 system.
Following preprocessing, the data set consists of 992 observational samples and 13,534 interventional samples over eight gene deletions.  As in~\cite{D16,greedy_SP_interventions}, we focused on learning the DAG structure on $24$ genes  that are transcription factors known to regulate expression of a variety of different genes, including one another \cite{G12}.

\begin{figure}[t!]
\centering
\subfigure[Effects of gene deletions]{\includegraphics[width=0.32\textwidth]{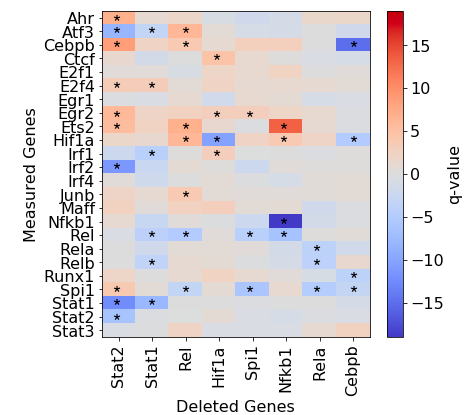}}
\subfigure[Recovery of gene deletion effects]{\includegraphics[width=0.32\textwidth]{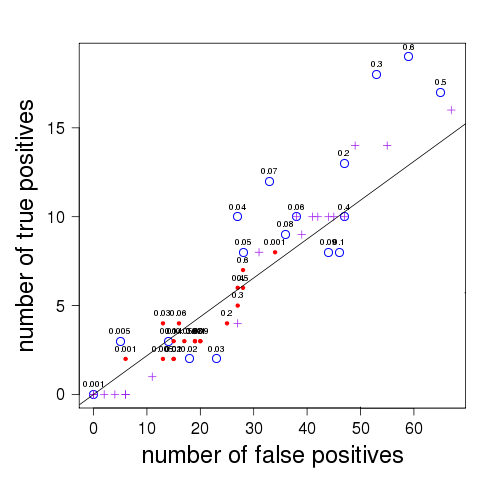}}
\subfigure[Network from stability selection]{\includegraphics[width=0.32\textwidth]{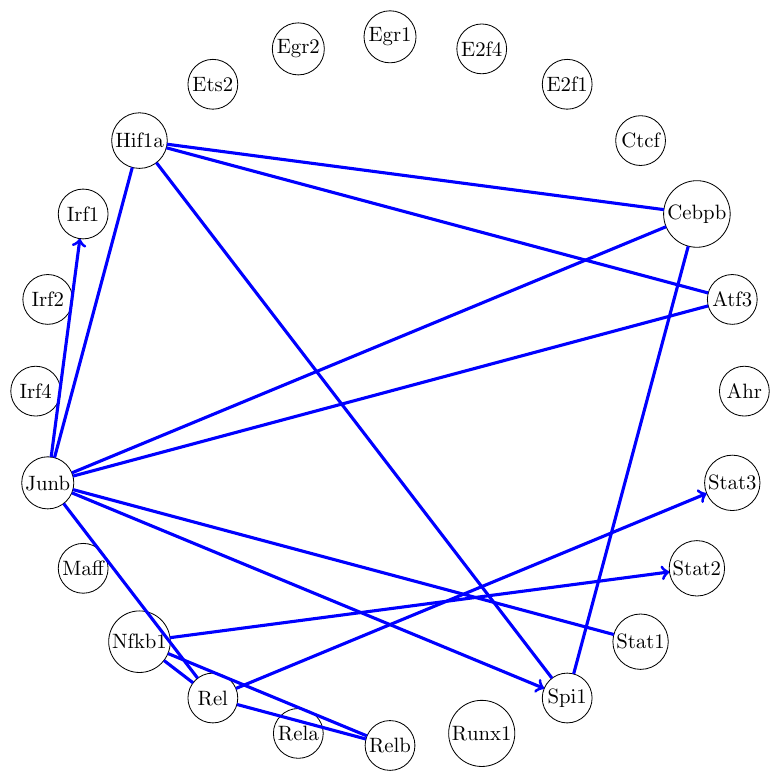}}
\caption{(a) Heatmap indicating the effect of each gene deletion on each measured gene;  
q-values with magnitude at least $3$ are marked with ``*''. (b) Performance of the causal network learned by Algorithm~\ref{alg: greedy sp depth and start control} with $d = 4$ and $r = 20$ (circle) as compared to the PC-algorithm (dot) and GES (cross) in predicting the effect of each intervention; line corresponds to random guessing. 
(c) The PDAG discovered by Algorithm~\ref{alg: greedy sp depth and start control} via stability selection using cutoff $0.6$.}

	\label{fig: Dixit}
\end{figure}

We used the observational samples to infer the DAG and the interventional samples to evaluate it. In particular, using the interventional data corresponding to a deletion of gene A we identified the genes that are downstream of gene A by testing whether the interventional distribution is significantly different from the observational distribution. For this, we used a Wilcoxon Rank-Sum test with p-value $\alpha = 0.05$, corresponding to a magnitude of at least $3$ in the q-value heat map depicted in Figure~\ref{fig: Dixit}(a); a positive q-value indicates that the gene expression level is increased by the gene deletion, whereas a negative value means that it is decreased. The accuracy of an estimated causal network is evaluated based on the edges adjacent to intervened nodes: an arrow from gene $A$ to gene $B$ in the learned network is considered a true positive if the expression of gene $B$ in the interventional distribution when targeting gene A is significantly different from the observational distribution,~e.g., there is a star in the $(A,B)$-entry in Figure~\ref{fig: Dixit}(a), and it is considered a false positive otherwise.
Using this metric, Figure~\ref{fig: Dixit}(b) compares the performance of Algorithm~\ref{alg: greedy sp depth and start control} with $d=4$ and $r = 20$, the PC-algorithm, and greedy equivalence search.  
Each point in Figure~\ref{fig: Dixit}(b) corresponds to the number of true positives and false positives in a DAG on the $24$ genes learned from the observational data with a fixed parameter setting.  
The fixed parameter is the significance level of the conditional independence test for Algorithm~\ref{alg: greedy sp depth and start control} and the PC-algorithm, and it is the $\ell_0$-penalization constant $c$ in the penalty $\lambda_n = c\log(n)$ used in the score function for greedy equivalence search. While the PC-algorithm performs similar to random guessing, the other two algorithms perform better, with the greedy sparsest permutation algorithm, i.e., Algorithm~\ref{alg: greedy sp depth and start control}, generally outperforming greedy equivalence search.

To get a sense of the corresponding gene regulatory network, in Figure~\ref{fig: Dixit}(c) we plotted the network constructed from our algorithm using stability selection~\cite{MB10}. We determined the cutoff parameter for stability selection by varying it between $0.5$ to $0.8$ and found that the resulting network was very robust in the range $[0.6, 0.7]$. The network shown in Figure~\ref{fig: Dixit}(c) corresponds to a cutoff of $0.6$, which is also within the recommended range given in~\cite{MB10}. 
In the supplementary material, we provide ROC plots using a q-value cutoff of $1$ instead of $3$, showing that our results and conclusions regarding the comparison of the different algorithms are robust with respect to the selection of q-value cutoff.

\section{Discussion}
\label{sec_discussion}
The greedy sparsest permutation algorithm, i.e., Algorithm~\ref{alg: greedy sp depth and start control}, with parameter choices $d = \infty$ and $r = 1$ is Algorithm~\ref{alg_triangle_SP}, which was shown to be consistent under strictly weaker conditions than faithfulness.  
Algorithm~\ref{alg_triangle_SP} is an approximation of Algorithm~\ref{alg_edge_SP}, which is further consistent under strictly weaker conditions than Algorithm~\ref{alg_triangle_SP}.  
The fact that Algorithm~\ref{alg_triangle_SP} is consistent under strictly weaker conditions than faithfulness was observed by its performance on simulated data in Section~\ref{sec: simulations}.  
On the other hand, Algorithm~\ref{alg_edge_SP} was not simulated since we must produce the entire polytope $\mathcal{A}_p(\CC)$ so as to recover its edge graph.  
A complete characterization of the edges of $\mathcal{A}_p(\CC)$ would thus be of use, so that Algorithm~\ref{alg_edge_SP} can be implemented without computing $\mathcal{A}_p(\CC)$ and its entire edge graph.  
Such an implementation would likely recover the true Markov equivalence class more often than any of the algorithms in Section~\ref{sec: simulations}.  
Further perspectives on this could be gained via a characterization of all distributions satisfying Assumption~\ref{ass: edge SP} or Assumption~\ref{ass: triangle SP}.  

We expect the greedy permutation-based approaches developed in this paper to be useful in a variety of settings. 
For instance, extensions of Algorithm~\ref{alg: greedy sp depth and start control} to the setting where a mix of observational and interventional data is available were presented in~\cite{greedy_SP_interventions,soft_IGSP,UTIGSP}, and they were implemented using kernel-based conditional independence tests~\cite{Fukumizu,Tillman} which are better able to deal with non-linear structural equations and non-Gaussian noise.
Extensions of Algorithm~\ref{alg: greedy sp depth and start control} to the causally insufficient setting are also being developed~\cite{GSPo}. 
In addition, it would be interesting to extend Algorithm~\ref{alg: greedy sp depth and start control} so as to accommodate cyclic graphs.  

Since passage to a greedy permutation-based algorithm is often motivated by a need to efficiently search through a state space that is super-exponential in size, it would be interesting to compare the computational efficiency of the algorithms discussed in Section~\ref{sec: simulations}.  
Such studies could be conducted using the \verb+CausalDAG+~Python package available at \url{https://github.com/uhlerlab/causaldag}, which provides an efficient implementation of Algorithm~\ref{alg: greedy sp depth and start control}.  

\section*{Acknowledgement}
Liam Solus was partially supported by an NSF Mathematical Sciences Postdoctoral Research Fellowship (DMS - 1606407), the Wallenberg AI, Autonomous Systems and Software Program (WASP) funded by the Knut and Alice Wallenberg Foundation, and Starting Grant (Etableringsbidrag) No.~2019-05195 from The Swedish Research Council (Vetenskapsr\aa{}det).
Caroline Uhler was partially supported by NSF (DMS-1651995), ONR (N00014-17-1-2147 and N00014-18-1-2765), IBM, a Sloan Fellowship and a Simons Investigator Award.  
The authors are grateful to Robert Castelo for helpful discussions, as well as two reviewers and an associate editor at Biometrika for their quick and thoughtful comments.	


\begin{appendix}

\section{Background Material and an Example}
\label{app: greedy SP}

\subsection{Background}
\label{app: background}
Here, we provide some definitions from graph theory and causal inference that we will use in the coming proofs.  
Given a DAG $\G := ([p], A)$ with node set $[p]:=\{1,2,\ldots,p\}$ and arrow set $A$, we associate to the nodes of $\G$ a random vector $(X_1,\ldots,X_p)$ with a probability distribution $\p$.  
An arrow in $A$ is an ordered pair of nodes $(i,j)$ which we will often denote by $i\rightarrow j$.  
A {directed path} in $\G$ from node $i$ to node $j$ is a sequence of directed edges in $\G$ of the form $i\rightarrow i_1 \rightarrow i_2 \rightarrow \cdots \rightarrow j$. 
A {path} from $i$ to $j$ is a sequence of arrows between $i$ and $j$ that connect the two nodes without regard to direction.  
The {parents} of a node $i$ in $\G$ is the collection $\parents_G(i) := \{ k \in[p] \,:\, k\rightarrow i \in A\}$, and the {ancestors} of $i$, denoted $\ancestors_\G(i)$, is the collection of all nodes $k\in[p]$ for which there exists a directed path from $k$ to $i$ in $\G$. We do not include $i$ in $\ancestors_\G(i)$. 
The {descendants} of $i$, denoted $\descendants_\G(i)$, is the set of all nodes $k\in[p]$ for which there is a directed path from $i$ to $k$ in $\G$, and the {nondescendants} of $i$ is the collection of nodes $\nondescendants_\G(i) := [p]\backslash(\descendants_\G(i)\cup\{i\})$. 
When the DAG $\G$ is understood from context we write $\parents(i)$, $\ancestors(i)$, $\descendants(i)$, and $\nondescendants(i)$, for the parents, ancestors, descendants, and nondescendants of $i$ in $\G$, respectively.  
The analogous definitions and notation will also be used for any set $S\subset[p]$.  
If two nodes are connected by an arrow in $\G$ then we say they are {adjacent}.  
A triple of nodes $(i,j,k)$ is called {unshielded} if $i$ and $j$ are adjacent, $k$ and $j$ are adjacent, but $i$ and $k$ are not adjacent.  
An unshielded triple $(i,j,k)$ forms an {immorality} if it is of the form $i\rightarrow j \leftarrow k$.  
In any triple, shielded or not, with arrows $i\rightarrow j\leftarrow k$, the node $j$ is called a {collider}.  
Given disjoint subsets $A,B,C\subset[p]$ with $A\cap B=\emptyset$, we say that $A$ is {$d$-connected} to $B$ given $C$ if there exist nodes $i\in A$ and $j\in B$ for which there is a path between $i$ and $j$ such that every collider on the path is in $\ancestors(C)\cup C$ and no non-collider on the path is in $C$.  
If no such path exists, we say $A$ and $B$ are {$d$-separated} given $C$.

\begin{example}
\label{ex: not all edges correspond to covered arrow flips}
\begin{figure}[t!]
\centering
\includegraphics[width = 0.7\textwidth]{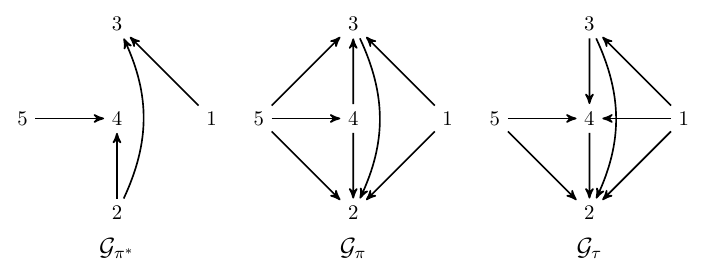}
\caption{An edge of a DAG associahedron that does not correspond to a covered edge flip.   The DAG associahedron $\mathcal{A}_p(\CC)$ is constructed for the conditional independence relations implied by the $d$-separation statements for $\G_{\pi^\ast}$ with $\pi^\ast = 15234$.  The DAGs $\G_\pi$ and $\G_\tau$ with $\pi = 15432$ and $\tau = 15342$ correspond to adjacent vertices in $\mathcal{A}_p(\CC)$, connected by the edge labeled by the transposition of $3$ and $4$.  The arrow between nodes $3$ and $4$ is not covered in either DAG $\G_\pi$ or $\G_\tau$.}  
\label{fig: not all edges correspond to covered arrow flips} 
\end{figure}
An example of a DAG associahedron containing an edge that does not correspond to a covered arrow reversal in either DAG labeling its endpoints can be constructed as follows:
Let $\G_{\pi^\ast}$ denote the left-most DAG depicted in Figure~\ref{fig: not all edges correspond to covered arrow flips}, and let $\CC$ denote those conditional independence relations implied by the $d$-separation statements for $\G_{\pi^\ast}$.  
Then for the permutations $\pi = 15432$ and $\tau = 15342$, the DAGs $\G_\pi$ and $\G_\tau$ label a pair of adjacent vertices of $\mathcal{A}_p(\CC)$ since $\pi$ and $\tau$ differ by the transposition of $3$ and $4$.  
This adjacent transposition corresponds to a reversal of the arrow between nodes $3$ and $4$ in $\G_\pi$ and $\G_\tau$.  
However, this arrow is not covered in either minimal independence map.  
We further note that this example shows that not all edges of $\mathcal{A}_p(\CC)$ can be described by covered arrow reversals even when $\CC$ is faithful to the sparsest minimal independence map, $\G_{\pi^\ast}$.  
\end{example}

\section{Proofs for Results on the Pointwise Consistency of the Greedy sparsest permutation algorithms}
\label{app: proofs of pointwise consistency guarantees for the greedy sparsest permutation algorithm}

\subsection{Proof of Lemma~\ref{thm: faithfulness}}
Suppose first that $\CC$ is not faithful to $\G$ and take any conditional independence statement $i\independent j \, \vert \, K$ that is not encoded by the $d$-separation statements in $\G$.  
Take $\pi$ to be any permutation in which $K\prec_\pi i \prec_\pi j \prec_\pi [p]\backslash(K\cup\{i,j\})$.  
Then $\G\not\leq\G_\pi$ since $i\independent j \, \vert \, K$ is encoded by the $d$-separations of $\G_\pi$ but not by the $d$-separations of $\G$.  

Conversely, suppose $\p$ is faithful to $\G$.  
By \cite[Theorem 9, Page 119]{Pearl_1988}, we know that $\p$ satisfies the Markov Assumption with respect to  $\G_\pi$  for any $\pi \in S_n$.  
So any conditional independence relation encoded by $\G_\pi$ holds for $\p$, which means it also holds for $\G$.  
Thus, $\G\leq \G_\pi$.
\hfill$\square$

\bigskip

To prove that Algorithm~\ref{alg_triangle_SP} is consistent under faithfulness we require a number of lemmas pertaining to the steps of the Chickering algorithm.  
For the convenience of the reader, we recall the Chickering algorithm in Algorithm~\ref{alg: chickering}.  
{
	\begin{algorithm}[b!]
		\caption{APPLY-EDGE-OPERATION}
		\label{alg: chickering}
		\LinesNumbered
		\DontPrintSemicolon
		\SetAlgoLined
		\SetKwInOut{Input}{Input}
		\SetKwInOut{Output}{Output}
		\Input{DAGs $\G$ and $\HH$ where $\G\leq \HH$ and $\G\neq \HH$.}
		\Output{A DAG $\G^\prime$ satisfying $\G^\prime\leq \HH$ that is given by reversing an edge in $\G$ or adding an edge to $\G$.}
		\BlankLine 
		Set $\G^\prime := \G$.
		\;
		While $\G$ and $\HH$ contain a node $Y$ that is a sink in both DAGs and for which $\parents_\G(Y) = \parents_\HH(Y)$, remove $Y$ and all incident edges from both DAGs.
		\;
		Let $Y$ be any sink node in $\HH$. 
		\;
		If $Y$ has no children in $G$, then let $X$ be any parent of $Y$ in $\HH$ that is not a parent of $Y$ in $\G$.  
		Add the edge $X\rightarrow Y$ to $\G^\prime$ and return $\G^\prime$. 
		\;
		Let $D\in\descendants_\G(Y)$ denote the (unique) maximal element from $\descendants_\G(Y)$ within $\HH$.  
		Let $Z$ be any maximal child of $Y$ in $\G$ such that $D$ is a descendant of $Z$ in $\G$.  
		\;
		If $Y\rightarrow Z$ is covered in $\G$, reverse $Y\rightarrow Z$ in $\G^\prime$ and return  $\G^\prime$.
		\;
		If there exists a node $X$ that is a parent of $Y$ but not a parent of $Z$ in $\G$, then add $X\rightarrow Z$ to $\G^\prime$ and return $\G^\prime$.
		\;
		Let $X$ be any parent of $Z$ that is not a parent of $Y$.   
		Add $X\rightarrow Y$ to $\G^\prime$ and return $\G^\prime$.  
	\end{algorithm}
}

\begin{lemma}
\label{lem: one edge at a time}
Suppose $\G\leq \HH$ such that the Chickering algorithm has reached step $5$ and selected the arrow $Y\rightarrow Z$ in $\G$ to reverse.  
If $\,Y\rightarrow Z$ is not covered in $\G$, then there exists a Chickering sequence 
$$
\left(\G = \G^0, \G^1, \G^2, \ldots, \G^N \leq \HH\right)
$$
in which $\G^N$ is produced by the reversal of $Y\rightarrow Z$, and for all $\,i =1,2,\ldots, N-1$, the DAG $\G^i$ is produced by an arrow addition via step $7$ or $8$ with respect to the arrow $Y\rightarrow Z$.  
\end{lemma}

\begin{proof}
Until the arrow $Y\rightarrow Z$ is reversed, the set $\descendants_\G(Y)$ and the node choice $D\in \descendants_\G(Y)$ remain the same.  
This is because steps $7$ and $8$  only add parents to $Y$ or $Z$  that are already parents of $Y$ or $Z$, respectively.  
Thus, we can always choose the same $Y$ and $Z$ until $Y\rightarrow Z$ is covered.  
\end{proof}

For an independence map $\G\leq \HH$, the Chickering algorithm first deletes all sinks in $\G$ that have precisely the same parents in $\HH$, and repeats this process for the resulting graphs until there is no sink of this type anymore.  
This is the purpose of step $2$ of the algorithm.  
If the adjusted graph is $\widetilde{\G}$, the algorithm then selects a sink node in $\widetilde{\G}$, which, by construction, must have fewer parents than the same node in $\HH$ and/or some children.  
The algorithm then adds parents and reverses arrows until this node has exactly the same parents as the corresponding node in $\HH$.  
The following lemma shows that this can be accomplished one sink node at a time.  
The proof is clear from the statement of the algorithm.  
\begin{lemma}
\label{lem: one sink at a time}
Let $\G\leq \HH$.  
If $Y$ is a sink node selectable in step $3$ of the Chickering algorithm then we may always select $Y$ each time until it is deleted by step $2$.  
\end{lemma}

We would like to see how the sequence of graphs produced in Chickering's algorithm relates to the DAGs $\G_\pi$ for a set of conditional independence relations $\CC$.  
In particular, we would like to see that if $\G_\pi\leq \G_\tau$ for permutations $\pi,\tau\in S_p$, then there is a sequence of moves given by Chickering's algorithm that passes through a sequence of minimal independence maps taking us from $\G_\pi$ to $\G_\tau$.  
To do so, we require an additional lemma relating independence maps and minimal independence maps.  
To state this lemma we need to consider the two steps within Algorithm~\ref{alg: chickering} in which arrow additions occur.  
We now recall these two steps:

\begin{enumerate}[(i)]
	\item Suppose $Y$ is a sink node in $\G\leq \HH$.  
	If $Y$ is also a sink node in $\G$, then choose a parent $X$ of $Y$ in $\HH$ that is not a parent of $Y$ in $\G$, and add the arrow $X\rightarrow Y$ to $\HH$.  
	\item If $Y$ is not a sink node in $\G$, then there exists an arrow $Y\rightarrow Z$ in $\G$ that is oriented in the opposite direction in $\HH$.  
	If $Y\rightarrow Z$ is covered, the algorithm reverses it.  
	If $Y\rightarrow Z$ is not covered, there exists (in $\G$)  either 
		\begin{enumerate}[(a)]
			\item a parent $X$ of $Y$ that is not a parent of $Z$, 
			in which case, the algorithm adds the arrow $X\rightarrow Z$. 
			\item a parent $X$ of $Z$ that is not a parent of $Y$,  
			in which case, the algorithm adds the arrow $X\rightarrow Y$.  
		\end{enumerate}
\end{enumerate}
\begin{lemma}
\label{lem: edge removals for chickerings algorithm}
Let $\CC$ be a graphoid and $\G_\pi\leq \G_\tau$ with respect to $\CC$.    
Then the common sink nodes of $\G_\pi$ and $\G_\tau$ all have the same incoming arrows.  
In particular, the Chickering algorithm needs no instance of arrow additions $(i)$ to move from $\G_\pi$ to $\G_\tau$.  
\end{lemma}

\begin{proof}
Suppose on the contrary that there exists some sink node $Y$ in $\G_\pi$ and there is a parent node $X$ of $Y$ in $\G_\tau$ that is not a parent node of $Y$ in $\G_\pi$.  
Since $Y$ is a sink in both permutations, then there exists linear extensions $\hat{\pi}$ and $\hat{\tau}$ of the partial orders corresponding to $\G_\pi$ and $\G_\tau$ for which $Y= \hat{\pi}_p$ and $Y = \hat{\tau}_p$.
By \cite[Theorem 7.4]{MUWY16}, we know that $\G_\pi = \G_{\hat{\pi}}$ and $\G_\tau = \G_{\hat{\tau}}$.  
In particular, we know that 
$X\not\independent Y\,\vert\, [p]\backslash\{X,Y\}$ in $\G_{\hat{\tau}}$ and $X\independent Y\,\vert\, [p]\backslash\{X,Y\}$ in $\G_{\hat{\pi}}$.  
However, this is a contradiction, since both of these relations cannot simultaneously hold.  
\end{proof}

\subsection{Lemmata for the Proof of Proposition~\ref{thm: consistent under faithfulness}} 
To prove Proposition~\ref{thm: consistent under faithfulness} we must first prove a few lemmas.
Throughout the remainder of this section, we use the following notation:
Suppose that $\G\leq \HH$ for two DAGs $\G$ and $\HH$ and that
$$
C = (\G^0:=\G, \G^1, \G^2,\ldots, \G^N:=\HH)
$$
is a Chickering sequence from $\G$ to $\HH$.  
We let $\pi^i\in S_p$ denote a linear extension of $\G^i$ for all $i=0,1,\ldots,N$.  
For any DAG $\G$ let $\CI(\G)$ denote the collection of conditional independence relations encoded by the $d$-separation statements in $\G$.  
\begin{lemma}
\label{lem: subDAG}
Suppose that $\G_\tau$ is a minimal independence map of a graphoid $\CC$.  
Suppose also that $\G\approx\G_\tau$ and that $\G$ differs from $\G_\tau$ only by a covered arrow reversal.  
If $\pi$ is a linear extension of $\G$ then $\G_\pi$ is a subDAG of $\G$.  
\end{lemma}

\begin{proof}
Suppose that $\G$ is obtained from $\G_\tau$ by the reversal of the covered arrow $x\rightarrow y$ in $\G_\tau$.  
Without loss of generality, we assume that $\tau = SxyT$ and $\pi = SyxT$ for some disjoint words $S$ and $T$ whose letters are collectively in bijection with the elements in $[p]\backslash\{x,y\}$.  
So in $\G_\pi$, the arrows going from $S$ to $T$, $x$ to $T$, and $y$ to $T$ are all the same as in $\G_\tau$.  
However, the arrows going from $S$ to $x$ and $S$ to $y$ may be different.  
So, to prove that $\G_\pi$ is a subDAG of $\G$ we must show that for each letter $s$ in the word $S$
\begin{enumerate}[(1)]
	\item if $s\rightarrow x\notin\G_\tau$ then $s\rightarrow x\notin\G_\pi$, and
	\item if $s\rightarrow y\notin\G_\tau$ then $s\rightarrow y\notin\G_\pi$.
\end{enumerate}
To see this, notice that if $s\rightarrow x\notin\G_\tau$, then $s\rightarrow y\notin\G_\tau$ since $x\rightarrow y$ is covered in $\G_\tau$.  
Similarly, if $s\rightarrow y\notin\G_\tau$ then $s\rightarrow x\notin\G_\tau$.  
Thus, we know that $s\independent x\,\vert\, S\backslash s$ and $s\independent y\,\vert\, (S\backslash s)x$ are both in the collection $\CC$.  
It then follows from the semigraphoid property (2) given in Section~\ref{sec_greedy_SP} that $s\independent x\,\vert\, (S\backslash s)y$ and $s\independent y\,\vert\, S\backslash s$ are in $\CC$ as well.  
Therefore, $\G_\pi$ is a subDAG of $\G$.  
\end{proof}

\begin{lemma}
\label{lem: edge deletion}
Let $\CC$ be a graphoid and let 
$$
C = (\G^0:=\G_\pi, \G^1, \G^2,\ldots, \G^N:=\G_\tau)
$$
be a Chickering sequence from a minimal independence map $\G_\pi$ of $\CC$ to another $\G_\tau$.  
If, for some index $0\leq i < N$, $\G^i$ is obtained from $\G^{i+1}$ by deletion of an arrow $x\rightarrow y$ in $\G^{i+1}$ then $x\rightarrow y$ is not in $\G_{\pi^{i+1}}$.  
\end{lemma}

\begin{proof}
Let $\pi^{i+1} = SxTyR$ be a linear extension of $\G^{i+1}$ for some disjoint words $S$, $T$, and $R$ whose letters are collectively in bijection with the elements in $[p]\backslash\{x,y\}$.  
Since $\G_{\pi^\ast}\leq\G^i\leq\G^{i+1}$ then 
$$
\CC\supseteq\CI(\G_\pi)\supseteq\CI(\G^i)\supseteq\CI(\G^{i+1}).
$$
We claim that $x\independent y\,\vert\, ST\in\CI(\G^i)\subseteq\CC$.  
Therefore, $x\rightarrow y$ cannot be an arrow in $\G_{\pi^{i+1}}$.  

First, since $\G^i$ is obtained from $\G^{i+1}$ by deleting the arrow $x\rightarrow y$, then $\pi^{i+1}$ is also a linear extension of $\G^i$.  
Notice, there is no directed path from $y$ to $x$ in $\G^i$, and so it follows that $x$ and $y$ are $d$-separated in $\G^i$ by $\pa_{\G^i}(y)$.  
Therefore, $x\independent y\,\vert\, \pa_{\G^i}(y)\in\CI(\G^i)$.  
Notice also that $\pa_{\G^i}(y)\subset ST$ and any path in $\G^i$ between $x$ and $y$ lacking colliders uses only arrows in the subDAG of $\G^i$ induced by the vertices $S\cup T\cup\{x,y\} = [p]\backslash R$.  
Therefore, $x\independent y\,\vert\, ST\in\CI(\G^i)$ as well.  
It follows that $x\independent y\,\vert\, ST\in\CC$, and so, by definition, $x\rightarrow y$ is not an arrow of $\G_{\pi^{i+1}}$.  
\end{proof}

\begin{lemma}
\label{lem: proper subDAG}
Suppose that $\CC$ is a graphoid and $\G_\pi$ is a minimal independence map with respect to $\CC$.  
Let 
$$
C = (\G^0:=\G_\pi, \G^1, \G^2,\ldots, \G^N:=\G_\tau)
$$
be a Chickering sequence from $\G_\pi$ to another minimal independence map $\G_\tau$ with respect to $\CC$.  
Let $i$ be the largest index such that $\G^i$ is produced from $\G^{i+1}$ by deletion of an arrow, and suppose that for all $i+1< k\leq N$ we have $\G_{\pi^k}=\G^k$.  
Then $\G_{\pi^{i+1}}$ is a proper subDAG of $\G^{i+1}$.  
\end{lemma}

\begin{proof}
By Lemma~\ref{lem: subDAG}, we know that $\G_{\pi^{i+1}}$ is a subDAG of $\G^{i+1}$.  
This is because $\pi^{i+1}$ is a linear extension of $\G^{i+1}$ and $\G^{i+1}\approx\G^{i+2} = \G_{\pi^{i+2}}$ and $\G^{i+1}$ differs from $\G^{i+2}$ only by a covered arrow reversal.  
By Lemma~\ref{lem: edge deletion}, we know that the arrow deleted in $\G^{i+1}$ to obtain $\G^i$ is not in $\G_{\pi^{i+1}}$.  
Therefore, $\G_{\pi^{i+1}}$ is a proper subDAG of $\G$.  
\end{proof}

Using these lemmas, we can now give a proof of Proposition~\ref{thm: consistent under faithfulness}.

\subsection{Proof of Proposition~\ref{thm: consistent under faithfulness}}
To see that (a) holds, notice since $\G_\pi\approx\G_\tau$ then by the transformational characterization of Markov equivalence given in \cite[Theorem 2]{C95}, we know there exists a Chickering sequence
$$
C := (\G^0:=\G_\pi, \G^1, \G^2,\ldots, \G^N:=\G_\tau)
$$
for which $\G^0\approx\G^1\approx\cdots\approx\G^N$ and $\G^i$ is obtained from $\G^{i+1}$ by the reversal of a covered arrow in $\G^{i+1}$ for all $0\leq i< N$.  
Furthermore, since $\G_\pi$ is class-s-minimal, and by Lemma~\ref{lem: subDAG}, we know that for all $0\leq i\leq N$
$$
\overline{\G}^i\supseteq\overline{\G}_{\pi^{i}}\supseteq\overline{\G}_\pi.  
$$
However, since $\G^i\approx\G_{\pi}$ and $\G_{\pi^{i}}$ is a subDAG of $\G^i$, then $\G^i= \G_{\pi^i}$ for all $i$.  
Thus, the desired weakly decreasing edgewalk along $\mathcal{A}_p(\CC)$ is 
$$
(\G_\pi = \G_{\pi^0},\G_{\pi^1},\ldots,\G_{\pi^{N-1}},\G_{\pi^N}=\G_\tau).
$$

To see that (b) holds, suppose that $\G_\pi\leq\G_\tau$ but $\G_\pi\not\approx\G_\tau$.  
Since $\G_\pi\leq\G_\tau$, there exists a Chickering sequence from $\G_\pi$ to $\G_\tau$ that uses at least one arrow addition.  
By Lemmas~\ref{lem: one edge at a time} and~\ref{lem: one sink at a time} we can choose this Chickering sequence such that it resolves one sink at a time and, respectively, reverses one covered arrow at a time. 
We denote this Chickering sequence by 
$$
C := (\G^0:=\G_\pi, \G^1, \G^2,\ldots, \G^N:=\G_\tau).  
$$
Let $i$ denote the largest index for which $\G^i$ is obtained from $\G^{i+1}$ by deletion of an arrow.  
Then by our choice of Chickering sequence we know that $\G^k$ is obtained from $\G^{k+1}$ by a covered arrow reversal for all $i<k<N$.  
Moreover, $\pi^i = \pi^{i+1}$, and so $\G_{\pi^i} = \G_{\pi^{i+1}}$.  
Furthermore, by Lemma~\ref{lem: subDAG} we know that $\G_{\pi^k}$ is a subDAG of $\G^k$ for all $i<k\leq N$.  

Suppose now that there exists some index $i+1<k<N$ such that $\G_{\pi^k}$ is a proper subDAG of $\G^k$.   
Without loss of generality, we pick the largest such index.  
It follows that for all indices $k<\ell\leq N$, $\G_{\pi^\ell} = \G^\ell$ and that
$$
\G^{k+1}\approx\G^{k+2}\approx\cdots\approx\G^N = \G_\tau.
$$
Thus, by \cite[Theorem 2]{C95}, there exists a weakly decreasing edgewalk from $\G_\tau$ to $\G^{k+1}$ on $\mathcal{A}_p(\CC)$.  
Since we chose the index $k$ maximally then $\G^k$ is obtained from $\G^{k+1}$ by a covered arrow reversal.  
Therefore, $\G_{\pi^k}$ and $\G_{\pi^{k+1}}$ are connected by an edge of $\mathcal{A}_p(\CC)$ indexed by a covered arrow reversal.  
Since $|\G^k| = |\G^{k+1}| = |\G_{\pi^{k+1}}|$ and $\G_{\pi^k}$ is a proper subDAG of $\G^k$, then the result follows.  

On the other hand, suppose that for all indices $i+1<k\leq N$, we have $\G_{\pi^k} = \G^k$.  Then this is precisely the conditions of Lemma~\ref{lem: proper subDAG}, and so it follows that $\G_{\pi^{i+1}}$ is a proper subDAG of $\G^{i+1}$.  
Since $\G^{i+1}$ is obtained from $\G^{i+2}$ by a covered arrow reversal, the result follows.  
\hfill$\square$

\subsection{Proof of Theorem~\ref{thm: alg_bic}} 
The proof is composed of two parts.  
We first prove that for any permutation $\pi$, in the limit of large $n$, $\hat{\G}_\pi$ is a minimal independence map of $\G_\pi$.  
We prove this by contradiction.   
Suppose $\hat{\G}_\pi \neq \G_\pi$.  
Since the Bayesian information criterion is a consistent scoring function \cite{H88}, in the limit of large $n$, $\hat{\G}_\pi$ is an independence map of the distribution. 
Since $\hat{\G}_\pi$ and $\G_\pi$ share the same permutation and $\G_\pi$ is a minimal independence map, then $\G_\pi \subset \hat{\G}_\pi$. 
Suppose now that there exists $(i,j) \in \hat{\G}_\pi$ such that $(i, j) \not\in \G_\pi$.  
Since $\G_\pi$ is a minimal independence map, we obtain that
$
i \independent j \, \vert \, {\pa}_{\G_\pi}(j).  
$
In Lemma~$7$ of \cite{C02}, it is shown that Bayesian scoring is locally consistent, and it follows from the first sentence of the proof therein that the Bayesian information criterion is also locally consistent.
Since the Bayesian information criterion is locally consistent, it follows that $\textrm{BIC}(\G_\pi, \hat{X}) > \textrm{BIC}(\hat{\G}_\pi, \hat{X})$.

Now we prove that for any two permutations $\tau$ and $\pi$ where $\G_\tau$ is connected to $\G_\pi$ by precisely one covered arrow reversal, in the limit of large~$n$, 
\begin{align*}
\textrm{BIC}(\G_\tau; \hat{X}) > \textrm{BIC}(\G_\pi; \hat{X}) \Leftrightarrow \vert \G_\tau \vert < \vert \G_\pi \vert,
\end{align*}
and
\begin{align*}
\textrm{BIC}(\G_\tau; \hat{X}) = \textrm{BIC}(\G_\pi; \hat{X}) \Leftrightarrow \vert \G_\tau \vert = \vert \G_\pi \vert.
\end{align*}
It suffices to prove 
\begin{align}
\vert \G_\tau \vert = \vert \G_\pi \vert \Rightarrow \textrm{BIC}(\G_\tau; \hat{X}) = \textrm{BIC}(\G_\pi; \hat{X}) \label{eq:score1}
\end{align}
and 
\begin{align}
\vert \G_\tau \vert < \vert \G_\pi \vert \Rightarrow \textrm{BIC}(\G_\tau; \hat{X}) > \textrm{BIC}(\G_\pi; \hat{X}). \label{eq:score2}
\end{align}
Eq.~\ref{eq:score1} is easily seen to be true using \cite[Theorem 2]{C95} as $\G_\pi$ and $\G_\tau$ are equivalent. 
For Eq.~\ref{eq:score2}, by Theorem~\ref{thm: meek's conjecture}, since $\G_\tau \leq \G_\pi$ there exists a Chickering sequence from $\G_\tau$ to $\G_\pi$ with at least one edge addition and several covered arrow reversals. 
For the covered arrow reversals, the Bayesian information criterion remains the same since the involved DAGs are equivalent.  
For the edge additions, the score necessarily decreases in the limit of large $n$ due to the increase in the number of parameters.   
This follows from the consistency of the Bayesian information criterion and the fact that DAGs before and after edge additions are both independence maps of $\mathbb{P}$.
In this case, the path taken in the triangle sparsest permutation algorithm using the Bayesian information criterion is the same as in the original triangle sparsest permutation algorithm.  
Since the triangle sparsest permutation algorithm is consistent, it follows that the triangle sparsest permutation algorithm with the Bayesian information criterion is also consistent.
\hfill$\square$ 

\subsection{Proof of Theorem~\ref{thm: faithfulness, SMR, and greedy SP assumptions}}
\begin{figure}[t!]
	\centering
	\includegraphics[width = 0.2\textwidth]{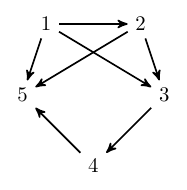}
	\caption{A sparsest DAG w.r.t.~the conditional independence relations $\CC$ given in the proof of Theorem~\ref{thm: faithfulness, SMR, and greedy SP assumptions}.}
	\label{fig: TSP weaker than faithfulness}
\end{figure}

It is quick to see that 
\begin{equation*}
\begin{split}
\textrm{faithfulness} \quad&\Longrightarrow \quad\textrm{triangle assumption} \\
\textrm{triangle assumption}\quad&\Longrightarrow\quad \textrm{edge assumption}	\\
 \textrm{edge assumption}\quad&\Longrightarrow\quad \textrm{sparsest Markov representation assumption}.	\\
\end{split}
\end{equation*}
The first implication is given by Theorem~\ref{cor: consistent under faithfulness}, and the latter three are immediate consequences of the definitions of the triangle, edge, and sparsest Markov representation assumptions.   
Namely, the triangle, edge, and sparsest Markov representation assumptions are each defined to be precisely the condition in which Algorithm~\ref{alg_triangle_SP}, Algorithm~\ref{alg_edge_SP}, and the sparsest permutation algorithm are, respectively, consistent.  
The implications then follow since each of the algorithms is a refined version of the preceding one in this order.  
Hence, we only need to show the strict implications. For each statement we identify a collection of conditional independence relations satisfying the former identifiability assumption but not the latter.  
For the first implication consider the collection of conditional independence relations
\begin{equation*}
\begin{split}
\CC = \{ 1\independent 5 \, \vert \, \{2,3\}, 
\quad
2\independent 4 \, \vert \, \{1,3\}, 
\quad
3\independent 5 \, \vert \, \{1,2,4\}, \\
\quad
1\independent 4 \, \vert \, \{2,3,5\}, 
\quad
1\independent 4 \, \vert \, \{2,3\}\}.
\end{split}
\end{equation*}
\noindent The sparsest DAG $\G_{\pi^*}$ with respect to $\CC$ is shown in Figure~\ref{fig: TSP weaker than faithfulness}. To see that $\CC$ satisfies the triangle assumption with respect to $\G_{\pi^*}$, we can use computer evaluation.  
To see that it is not faithful with respect to $\G_\pi^\ast$, notice that $1\independent 5 \, \vert \, \{2,3\}$ and $1\independent 4 \, \vert \, \{2,3,5\}$ are both in $\CC$, but they are not implied by $\G_\pi^\ast$.  
We also remark that $\CC$ is not a semigraphoid since the semigraphoid property (2) given in Section~\ref{sec_greedy_SP} applied to the conditional independence relations $1\independent 5 \, \vert \, \{2,3\}$ and $1\independent 4 \, \vert \, \{2,3,5\}$ implies that $1\independent 5 \, \vert \, \{2,3,4\}$ should be in $\CC$.  

For the second implication consider the collection of conditional independence relations
$$
\mathcal{D} = \{ 1\independent 2 \, \vert \, \{4\}, 
\quad
1\independent 3 \, \vert \, \{2\}, 
\quad
2\independent 4 \, \vert \, \{1,3\}\}
$$
and initialize Algorithm~\ref{alg_triangle_SP} at the permutation $\pi := 1423$.  
	\begin{figure}[b]
	\centering
	\includegraphics[width = 0.75\textwidth]{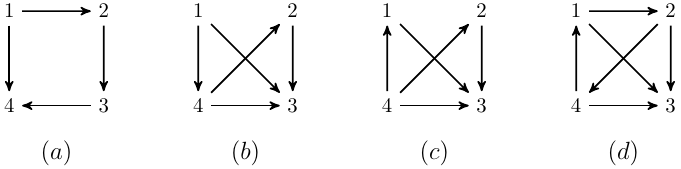}
	\caption{The four minimal independence maps with respect to the conditional independence relations $\mathcal{D}$ described in the proof of Theorem~\ref{thm: faithfulness, SMR, and greedy SP assumptions}.}
	\label{fig: ESP weaker than TSP}
	\end{figure}
A sparsest DAG $\G_{\pi^*}$ with respect to $\mathcal{D}$ is given in Figure~\ref{fig: ESP weaker than TSP}(a), and the initial minimal independence map $\G_\pi$ is depicted in Figure~\ref{fig: ESP weaker than TSP}(b).  
Notice that the only covered arrow in $\G_\pi$ is $1 \rightarrow 4$ , and reversing this covered arrow produces the permutation $\tau = 4123$; the corresponding DAG $\G_\tau$ is shown in Figure~\ref{fig: ESP weaker than TSP}(c).  
The only covered arrows in $\G_\tau$ are $4\rightarrow 1$ and $4\rightarrow 2$.  
Reversing $4\rightarrow 1$ returns us to $\G_\pi$, which we already visited, and reversing $4 \rightarrow 2$ produces the permutation $\sigma = 2143$; the associated DAG $\G_\sigma$ is depicted in Figure~\ref{fig: ESP weaker than TSP}(d).  
Since the only DAGs connected to $\G_\pi$ and $\G_\tau$ via covered arrow flips have at least as many edges as $\G_\pi$ and $\G_\tau$, then Algorithm~\ref{alg_triangle_SP} is inconsistent, and so the triangle assumption does not hold for $\CC$.  
On the other hand, we can verify computationally that Algorithm~\ref{alg_edge_SP} is consistent with respect to $\mathcal{D}$, meaning that the edge assumption holds.

	\begin{figure}[t]
	\centering
	\includegraphics[height = 1.5in]{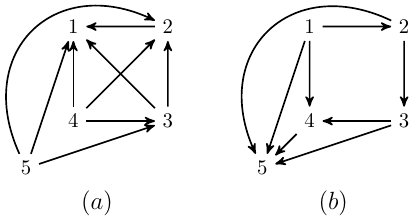}
	\caption{The initial minimal independence map and the sparsest minimal independence map with respect to the conditional independence relations $\mathcal{E}$ described in the proof of Theorem~\ref{thm: faithfulness, SMR, and greedy SP assumptions}.}
	\label{fig: SMR weaker than ESP}
	\end{figure}

Finally, for the last implication consider the collection of conditional independence relations
$$
\mathcal{E} = \{ 1\independent 3 \, \vert \, \{2\}, 
\quad
2\independent 4 \, \vert \, \{1,3\}, 
\quad
4\independent 5\},
$$
and the initial permutation $\pi = 54321$.  The initial DAG $\G_\pi$ and a  sparsest DAG $\G_{\pi^\ast}$  are depicted in Figures~\ref{fig: SMR weaker than ESP}(a) and (b), respectively.   
It is not hard to check that any DAG $\G_\tau$ that is edge adjacent to $\G_{\pi}$ is a complete graph.  
Thus, the sparsest Markov representation assumption holds for $\mathcal{E}$ but not the edge assumption.  
\hfill$\square$

\begin{figure}[b!]
	\centering
	\includegraphics[width = 0.2\textwidth]{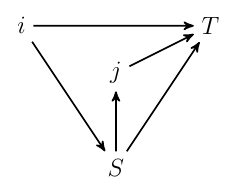}
	\caption{This diagram depicts the possible arrows between the node sets $\{i\}, \{j\}, S,$ and $T$ for the minimal independence map $\G_\pi$ considered in the proof of Theorem~\ref{thm: TSP and adjacency faithfulness}.}
	\label{fig: TSP and adjacency faithfulness diagram}
\end{figure}

\subsection{Proof of Theorem~\ref{thm: TSP and adjacency faithfulness}}
Let $\p$ be a semigraphoid, and let $\CC$ denote the conditional independence relations entailed by $\p$.  
Suppose for the sake of contradiction that Algorithm~\ref{alg_triangle_SP} is consistent with respect to $\CC$, but $\p$ fails to satisfy adjacency faithfulness with respect to a sparsest DAG $\G_\pi^\ast$.  
Then there exists some conditional independence relation $i\independent j\,\vert\, S$ in $\CC$ such that $i\rightarrow j$ is an arrow of $\G_\pi^\ast$.  
Now let $\pi$ be any permutation respecting the concatenated ordering $iSjT$ where $T = [p]\setminus(\{i,j\}\cup S)$.  
Then our goal is to show that any covered arrow reversal in $\G_\pi$ that results in a minimal independence map $\G_\tau$ with strictly fewer edges than $\G_\pi$ must satisfy the condition that $i\rightarrow j$ is not an arrow in $\G_\tau$.  

\begin{figure}[t!]
	\centering
	\includegraphics[width = 0.7\textwidth]{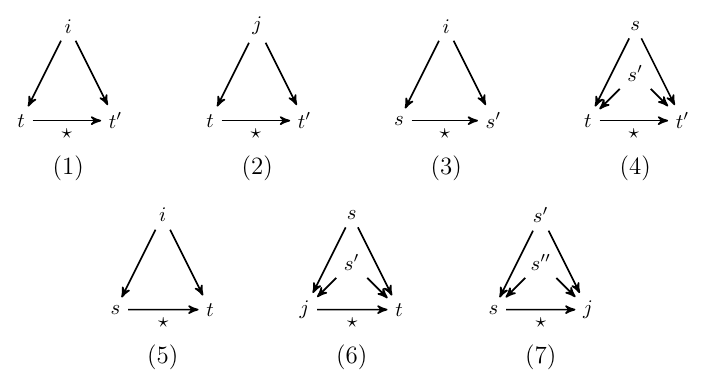}
	\caption{The possible non-trivially covered arrows between the node sets $\{i\}, \{j\}, S,$ and $T$ for the minimal independence map $\G_\pi$ considered in the proof of Theorem~\ref{thm: TSP and adjacency faithfulness} are labeled with the symbol $\star$.  Here, we take $s,s^\prime, s^{\prime\prime}\in S$ and $t,t^\prime\in T$.}
	\label{fig: TSP and adjacency faithfulness covered arrows}
	\end{figure}

First, we consider the possible types of covered arrows that may exist in $\G_\pi$.  
To list these, it will be helpful to look at the diagram depicted in Figure~\ref{fig: TSP and adjacency faithfulness diagram}. 
Notice first that we need not consider any trivially covered arrows, since such edge reversals do not decrease the number of arrows in the minimal independence maps.  
Any edge $i\rightarrow S$ or $i\rightarrow T$ is trivially covered, so the possible cases of non-trivially covered arrows are exactly the covered arrows given in Figure~\ref{fig: TSP and adjacency faithfulness covered arrows}.  
In this figure, each covered arrow to be considered is labeled with the symbol $\star$. 
Notice that the claim is trivially true for cases (1) -- (4); i.e., any covered arrow reversal resulting in edge deletions produces a minimal independence map $\G_\tau$ for which $i \rightarrow j$ is not an arrow of $\G_\tau$.  

Case (5) is also easy to see.  
Recall that $\pi = is_1\cdots s_kjt_1\cdots t_m$ where $S := \{s_1,\ldots, s_k\}$ and $T := \{t_1,\ldots, t_k\}$, and that reversing the covered arrow in case (5) results in an edge deletion.  
Since $s\rightarrow t$ is covered, then there exists a linear extension $\tau$ of $\G_\pi$ such that $s$ and $t$ are adjacent in $\tau$.  
Thus, either $j$ precedes both $s$ and $t$ or $j$ follows both $s$ and $t$ in $\tau$.  
Recall also that by \cite[Theorem 7.4]{MUWY16} we known $\G_\tau = \G_\pi$.  
Thus, reversing the covered arrow $s\rightarrow t$ in $\G_\tau = \G_\pi$ does not add in $i\rightarrow j$.  

To see the claim also holds for cases (6) and (7), we utilize the semigraphoid property (2) given in Section~\ref{sec_greedy_SP}.  
It suffices to prove the claim for case (6).  
So suppose that reversing the $\star$-labeled edge $j\rightarrow t$ from case (6) results in a minimal independence map with fewer arrows.  
We simply want to see that $i\rightarrow j$ is still a non-arrow in this new DAG.  
Assuming once more that $\pi = is_1\cdots s_kjt_1\cdots t_m$, by \cite[Theorem 7.4]{MUWY16} we can, without loss of generality, pick $t := t_1$.  
Thus, since $i\independent j \,\vert\, S$ and $j\rightarrow t$ is covered, then $i\independent t \,\vert\, S\cup\{j\}$.  
By the semigraphoid property (2), we then know that $i\independent j\,\vert\, S\cup\{t\}$.  
Thus, the covered arrow reversal $j\leftarrow t$ produces a permutation $\tau = is_1\cdots s_kt_1jt_2\cdots t_m$, and so $i\rightarrow j$ is not an arrow in $\G_\tau$.  
This completes all cases of the proof. 

	\begin{figure}[b]
	\centering
	\includegraphics[width = 0.3\textwidth]{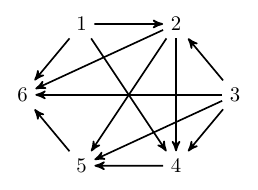}
	\caption{A sparsest DAG for the conditional independence relations $\CC$ considered in the proof of Theorem~\ref{thm: TSP and adjacency faithfulness}.}
	\label{fig: orientation faithfulness sparsest DAG}
\end{figure}

To complete the proof, we provide an example of a distribution $\mathbb{P}$ that satisfies the triangle assumption but not orientation faithfulness. Consider any probability distribution entailing the conditional independence relations
$$
\CC = \{ 1\independent 3, 
\quad
1\independent 5 \, \vert \, \{2,3,4\}, 
\quad
4\independent 6 \, \vert \, \{1,2,3,5\}, 
\quad
1\independent 3 \, \vert \, \{2,4,5,6\}
\}.
$$
For example, $\CC$ can be faithfully realized by a regular Gaussian. 
From left-to-right, we label these conditional independence relations as $c_1,c_2,c_3,c_4$.  
For the collection $\CC$, a sparsest DAG $\G_{\pi^\ast}$ is depicted in Figure~\ref{fig: orientation faithfulness sparsest DAG}. Note that since there is no equally sparse or sparser DAG that is Markov with respect to $\p$ then $\p$ satisfies the sparsest Markov representation assumption with respect to $\G_{\pi^\ast}$. Notice also that the conditional independence relation $c_4$ does not satisfy the orientation faithfulness assumption with respect to $\G_{\pi^\ast}$.
Moreover, if $\G_\pi$ entails $c_4$, then the subDAG on the nodes $\pi_1,\ldots, \pi_5$ forms a complete graph.  
Thus, by \cite[Theorem 2]{C95}, we can find a sequence of covered arrow reversals preserving edge count such that after all covered arrow reversals, $\pi_5 = 6$.  
Then transposing the entries $\pi_5\pi_6$ produces a permutation $\tau$ in which $c_3$ holds.  
Therefore, the number of arrows in $\G_\pi$ is at least the number of arrows in $\G_\tau$.  
Even more, $\G_\tau$ is an independence map of $\G_{\pi^\ast}$, i.e., $\G_{\pi^\ast}\leq \G_\tau$.  
So by Proposition~\ref{thm: consistent under faithfulness}, there exists a weakly decreasing edge walk determined by covered arrow reversals along $\mathcal{A}_p(\CC)$ taking us from $\G_\tau$ to $\G_\pi^\ast$.  
Thus, we conclude that $\p$ satisfies the triangle assumption, but not orientation faithfulness.  
\hfill$\square$


\section{Proofs for Results on the Uniform Consistency of the Greedy sparsest permutation algorithm}
\label{app: proofs for results on uniform consistency}

\subsection{Lemmata for the Proof of Theorem~\ref{thm_high_dim_consist}}
\label{subsec: high dim lemmata}
To prove Theorem~\ref{thm_high_dim_consist}, we require a pair of lemmas, the first of which shows that the conditioning sets in the triangle sparsest permutation algorithm can be restricted to parent sets of covered arrows.
\begin{lemma} 
\label{lm:flip}
Suppose that the data-generating distribution $\PP$ is faithful to $\mathcal{G}^*$. 
Then for any permutation $\pi$ and any covered arrow $i\rightarrow j$ in $\mathcal{G}_\pi$ it holds that
\begin{enumerate}
\item[(a)] $i \independent k \vert (S' \union \lbrace j \rbrace) \setminus \lbrace k \rbrace$ \;if and only if\;  $i \independent k \vert (S \union \lbrace j \rbrace) \setminus \lbrace k \rbrace$,
\item[(b)] $j \independent k \vert S' \setminus \lbrace k \rbrace$ \;if and only if\; $j \independent k \vert S \setminus \lbrace k \rbrace$,
\end{enumerate}
for all $k\in S$, where $S$ is the set of common parent nodes of $i$ and $j$, and $S' = \lbrace a : a <_\pi \max_\pi(i, j) \rbrace$.
\end{lemma}

\begin{proof}
Let $\pa_{\mathcal{G}_\pi}(j)$ be the set of parent nodes of node $j$ in the DAG $\mathcal{G}_\pi$. 
Let $k \in S$ and let $\PP_1$ denote the joint distribution of $(X_i, X_j, X_k)$ conditioned on $S\setminus\{k\}$ and $\PP_2$ the joint distribution of $(X_i, X_j, X_k)$ conditioned on $S'$. Then the claimed statements boil down to 
\begin{enumerate}
\item[(a)] $j \independent k\; \text{under distribution}\; \PP_1$ $\Leftrightarrow$ $j \independent k\; \text{under distribution}\; \PP_2$;
\item[(b)] $ i \independent k \vert j\; \text{under distribution}\; \PP_1$ $\Leftrightarrow$  $i \independent k \vert j\; \text{under distribution}\; \PP_2$.
\end{enumerate}

Note that
\begin{align*}
\begin{split}
\PP_1(X_i, X_j, X_k) & := \PP(X_i, X_j, X_k \vert X_{S \setminus \lbrace k \rbrace}) = \PP(X_i, X_j \vert X_{S})  \PP_1(X_k).
\end{split}
\end{align*}
Similarly, the Markov assumption of $\PP$ with respect to $G_\pi$ implies that
\begin{align*}
\begin{split}
\PP_2(X_i, X_j, X_k) &= \PP(X_i, X_j \vert X_{S'})  \PP_2(X_k) = \PP(X_i, X_j \vert X_{S})  \PP_2(X_k).
\end{split}
\end{align*}
Hence, $\PP_1(X_j \vert X_k) = \PP_2(X_j \vert X_k)$, $\PP_1(X_i \vert X_j, X_k) = \PP_2(X_i \vert X_j, X_k)$. This completes the proof since $X_a \independent X_b \vert X_C\; \text{under some distribution}\; \tilde{\PP}$ if and only if $\tilde{\PP}(X_a \vert X_b = z_1, X_C) = \tilde{\PP}(X_a \vert X_b = z_2, X_C)$ for all $z_1$ and $z_2$ in the sample space. 
\end{proof}

The second lemma we require was first proven in \cite[Lemma 3]{KB07} and is here restated for the sake of completeness.
\begin{lemma}
\label{lm:error}
\cite[Lemma 3]{KB07}
Suppose that assumption (4) holds, and let $z_{i,j \vert S}$ be the z-transform of the partial correlation coefficient $\rho_{i,j \vert S}$.  
Then
\begin{equation*}
\begin{split}
\PP&[\vert \hat{z}_{i,j \vert S} - z_{i,j \vert S} \vert > \gamma]\leq
\mathcal{O}(n - \vert S \vert)*\Phi, \mbox{ where}\\
\Phi&=\left[ \exp\left\{(n - 4 - \vert S \vert) \log \left(\frac{4 - (\gamma / L)^2}{4 + (\gamma / L)^2}\right)\right\} + \exp\{-C_2(n - \vert S \vert)\}\right],
\end{split}
\end{equation*}
where $C_2$ is some constant such that $0 < C_2 < \infty$ and 
$$
L = 1 / (1 - (1+M)^2 / 4),
$$
in which $M$ is defined such that it satisfies assumption (4).
\end{lemma}

Provided with Lemmas~\ref{lm:flip}~and~\ref{lm:error}, we can then prove Theorem~\ref{thm_high_dim_consist}.

\subsection{Proof of Theorem~\ref{thm_high_dim_consist}}
\label{subsubsec: proof of high dimensional consistency theorem}
For any initial permutation $\pi_0$, we let $L_{\pi_0}$ denote the set of tuples $(i, j, S)$ used for partial correlation testing in the estimation of the initial permtuation DAG $\mathcal{G}_{\pi_0}$.  
That is, 
\begin{align*}
L_{\pi_0} := \Big\lbrace (i, j, S): S = \lbrace k : \pi_0(k) \leq \max( \pi_0(i), \pi_0(j)) \rbrace \setminus \lbrace i,j \rbrace \Big\rbrace.  
\end{align*}
Given a DAG $\G$ and a node $i$ we let $\adj(\G, i)$ denote the collection of nodes that share an arrow with node $i$ in $G$.  
We then let $K_{\pi_0}$ denote the collection of tuples $(i, j, S)$ that will be used in the partial correlation testing done in step (2) of Algorithm~\ref{alg:grsp}; i.e. 
\begin{align*}
K_{\pi_0} := \underset{(i,j) \in \overline{\mathcal{G}}_{\pi_0}}{\bigcup} \Big\lbrace (k, l, S): k \in \lbrace i, j\rbrace,\; l \in \adj(\mathcal{G}_{\pi_0}, i) \intersection \adj(\mathcal{G}_{\pi_0}, j),\; S \neq \emptyset,\; \text{and}\; \\ 
S \subseteq \lbrace \adj(\mathcal{G}_{\pi_0}, i) \intersection \adj(\mathcal{G}_{\pi_0}, j) \rbrace \union \lbrace i, j \rbrace \Big\rbrace.
\end{align*}
It follows from Lemma~\ref{lm:flip}, that when flipping a covered edge $i\rightarrow j$ in a minimal independence map $\G_{\tilde{\pi}}$, it is sufficient to calculate the partial correlations $\rho_{a,b \vert C}$ where 
\begin{align*}
\begin{split}
(a, b, C) \in & \Big\lbrace (a, b, C) : a = i,\; b \in {\pa}_i(\mathcal{G}_{\tilde{\pi}}),\; C = {\pa}_i(\mathcal{G}_{\tilde{\pi}}) \union \lbrace j \rbrace \setminus \lbrace b \rbrace \Big\rbrace \union \\
& \Big\lbrace (a, b, C) : a = j,\; b \in {\pa}_i(\mathcal{G}_{\tilde{\pi}}),\; C = {\pa}_i(\mathcal{G}_{\tilde{\pi}}) \setminus \lbrace b \rbrace \Big\rbrace.  
\end{split}
\end{align*}
In particular, we have that $(a, b, C) \in K_{\tilde{\pi}}$.

Because of the skeletal inclusion $\overline{\mathcal{G}}_{\tilde{\pi}} \subseteq \overline{\mathcal{G}}_{\pi_0}$, it follows that $K_{\tilde{\pi}} \subseteq K_{\pi_0}$ and hence $(a, b, C) \in K_{\pi_0}$. In addition, for all partial correlations $\rho_{a,b \vert C}$ used for constructing the initial DAG $\mathcal{G}_{\pi_0}$, we know that $(a, b, C) \in L_{\pi_0}$.  
Therefore, for all partial correlations $(a, b, C)$ used in the algorithm, we have: 
\begin{align*}
(a, b, C) \in K_{\pi_0} \union L_{\pi_0}.
\end{align*}

Let $E_{i,j \vert S}$ be the event where an error occurs when doing partial correlation testing of $i \independent j \vert S$, and suppose that $\alpha$ is the significance level when testing this partial correlation.  
Then we see that $E_{i,j \vert S}$ corresponds to:
\begin{align*}
\begin{split}
(n - \vert S \vert - 3)^{1/2} \vert \hat{z}_{i,j \vert S} \vert > \Phi^{-1}(1 - \alpha / 2), \quad & \quad \text{when}\; z_{i,j \vert S} = 0; \\
(n - \vert S \vert - 3)^{1/2} \vert \hat{z}_{i,j \vert S} \vert \leq \Phi^{-1}(1 - \alpha / 2), \quad & \quad \text{when}\; z_{i,j \vert S} \neq 0.
\end{split}
\end{align*}
Choosing $\alpha_n = 2(1 - \Phi(n^{1 / 2} c_n / 2))$ it follows under assumption (4) that
\begin{align*}
\PP[E_{i,j \vert S}] \leq \PP[\vert \hat{z}_{i,j \vert S} - z_{i,j \vert S} \vert > (n / (n - \vert S \vert - 3))^{1/2} c_n / 2].
\end{align*}
Now, by (2) we have that $\vert S \vert \leq p = \mathcal{O}(n^a)$.  
Hence it follows that
\begin{align*}
\PP[E_{i,j \vert S}] \leq \PP[\vert \hat{z}_{i,j \vert S} - z_{i,j \vert S} \vert >  c_n / 2].
\end{align*}
Then, Lemma~\ref{lm:error} together with the fact that $\log(\frac{4 - \delta^2}{4 + \delta^2}) \sim - \delta^2 / 2$ as $\delta \rightarrow 0$, imply that
\begin{align}
\label{eq:error upper bound}
\PP[E_{i,j \vert S}] \leq \mathcal{O}(n - \vert S \vert) \exp\{-c' (n - \vert S \vert)c_n^2\} \leq \mathcal{O}\left( \exp(\log n - c n^{1 - 2\ell}) \right)
\end{align}
for some constants $c,c'>0$. Since the DAG estimated using Algorithm~\ref{alg:grsp} is not consistent when at least one of the partial correlation tests is not consistent, then the probability of inconsistency can be estimated as follows:
\begin{align} 
\label{eq:spbnd}
\begin{split}
\PP[\text{an error occurs in Algorithm~\ref{alg:grsp}}] & \leq \PP\left(\underset{i,j,S \in K_{\hat{\pi}} \union L_{\hat{\pi}}}{\bigcup} E_{i,j \vert S}\right) \\
& \leq \vert K_{\hat{\pi}} \union L_{\hat{\pi}} \vert \left(\underset{i,j,S \in K_{\hat{\pi}} \union L_{\hat{\pi}}}{\sup} \PP(E_{i,j \vert S})\right).
\end{split}
\end{align}
Next note that assumption (3) implies that the size of the set $\adj(\mathcal{G}_{\pi_0}, i) \union \adj(\mathcal{G}_{\pi_0}, j)$ is at most $d_{\pi_0}$.  
Therefore, $\vert K_{\pi_0} \vert \leq p^2 \cdot d_{\pi_0} \cdot 2^{d_{\pi_0}}$ and $\vert L_{\pi_0} \vert \leq p^2$.  
Thus, we see that
\begin{align*}
\vert K_{\hat{\pi}} \union L_{\hat{\pi}} \vert \leq \vert K_{\hat{\pi}} \vert + \vert L_{\hat{\pi}} \vert \leq (2^{d_{\pi_0}} \cdot d_{\pi_0} + 1) p^2.  
\end{align*}
Therefore, the left-hand-side of inequality~(\ref{eq:spbnd}) is upper-bounded by
\begin{align*}
(2^{d_{\pi_0}} \cdot d_{\pi_0} + 1) p^2 \left(\underset{i,j,S \in K_{\hat{\pi}} \union L_{\hat{\pi}}}{\sup} \PP(E_{i,j \vert S})\right).
\end{align*}
Combining this observation with the upper-bound computed in~(\ref{eq:error upper bound}), we obtain that the left-hand-side of~(\ref{eq:spbnd}) is upper-bounded by
\begin{align*}
(2^{d_{\pi_0}} \cdot d_{\pi_0} + 1) p^2 \mathcal{O}(\exp(\log n - c n^{1 - 2l})) \leq\\
\mathcal{O}(\exp(d_{\pi_0} \log 2 + 2 \log p + \log d_{\pi_0} + \log n - c n^{1 - 2\ell})).
\end{align*}
By assumptions (3) and (4) it follows that $n^{1 - 2\ell}$ dominates all terms in this bound.  
Thus, we conclude that
\begin{align*}
\PP[\text{estimated DAG is consistent}] \geq 1 - \mathcal{O}(\exp(- c n^{1 - 2\ell})).
\end{align*}
\hfill$\square$

The proof of Theorem~\ref{thm: NBMDA is equivalent to MDA in oracle setting} is based on the following lemma. 
\begin{lemma} 
\label{lm:mmd}
Let $\PP$ be a distribution on $[p]$ that is faithful to a DAG $\mathcal{G}$, and let $\PP_S$ denote the marginal distribution on $S\subset [p]$. 
Let $G_S$ be the undirected graphical model corresponding to $\PP_S$, i.e., the edge $\{i,j\}$ is in $G_S$ if and only if $\rho_{i,j \vert (S \setminus \lbrace i, j \rbrace)} \neq 0$.  
Then $G_{S \setminus \lbrace k \rbrace}$ can be obtained from $G_S$ as follows: 
\begin{enumerate}
	\item for all $i,j\in  \adj(G_S, k)$, if $\{i,j\}$ is not an edge in $G_S$, then add $\{i,j\}$.  Otherwise, $\{i,j\}$ is an edge of $G_{S \setminus \lbrace k \rbrace}$ if and only if $\vert \rho_{i,j \vert S \setminus \lbrace i, j, k \rbrace}\vert \neq 0$.
	\item for all $i,j\notin  \adj(G_S, k)$, $\{i,j\}$ is an edge of $G_{S \setminus \lbrace k \rbrace}$ if and only if $\{i,j\}$ is an edge in $G_S$.
\end{enumerate}
\end{lemma}

\begin{proof}
First, we prove:
\begin{align*}
\begin{split}
\textrm{For } i, j \not\in \adj(G_S, k):\quad & (i,j) \textrm{ is an edge in } G_{S \setminus \lbrace k \rbrace}\; \textrm{iff}\; (i,j) \textrm{ is an edge in } G_S. \\
\end{split}
\end{align*}
Suppose at least one of $i$ or $j$ are not adjacent to node $k$ in $G_S$.  
Without loss of generality, we assume $i$ is not adjacent to $k$ in $G_S$; this implies that $\rho_{i,k \vert S \setminus \lbrace i, k \rbrace} = 0$. 
To prove the desired result we must show that 
\begin{align*}
\rho_{i, j \vert S \setminus \lbrace i, j \rbrace} = 0 \Leftrightarrow \rho_{i, j \vert S \setminus \lbrace i, j, k \rbrace} = 0.
\end{align*}
To show this equivalence, first suppose that $\rho_{i, j \vert S \setminus \lbrace i, j \rbrace} = 0$ but $\rho_{i, j \vert S \setminus \lbrace i, j, k \rbrace} \neq 0$.  
This implies that there is a path $P$ between $i$ and $j$ through $k$ such that nodes $i$ and $j$ are d-connected given $S \setminus \lbrace i, j, k \rbrace$ and d-separated given $S \setminus \lbrace i, j\rbrace$.  
This implies that $k$ is a non-collider along $P$. Define $P_i$ as the path connecting $i$ and $k$ in the path $P$ and $P_j$ the path connecting $j$ and $k$ in $P$. Then the nodes $i$ and $j$ are d-connected to $k$ given $S \setminus \lbrace i, k \rbrace$ and $S \setminus \lbrace j, k \rbrace$ respectively, by using $P_i$ and $P_j$. Since $j$ is not on $P_i$, clearly $i$ and $k$ are also d-connected given $S \setminus \lbrace i, j, k \rbrace$ through $P_i$, and the same holds for $j$.

Conversely, suppose that 
$\rho_{i, j \vert S \setminus \lbrace i, j, k \rbrace} = 0$ but $\rho_{i, j \vert S \setminus \lbrace i, j \rbrace} \neq 0$.  
Then there exists a path $P$ that $d$-connects nodes $i$ and $j$ given $S \setminus \lbrace i, j \rbrace$, while $i$ and $j$ are d-separated given $S \setminus \lbrace i, j, k \rbrace$. 
Thus, one of the following must occur: 
\begin{enumerate}
\item $k$ is a collider on the path $P$, or 
\item Some node $\ell\in\text{an}(S \setminus \lbrace i, j \rbrace) \setminus \text{an}(S \setminus \lbrace i, j, k \rbrace)$ is a collider on $P$.
\end{enumerate}
For case (2), there must exist a path: $\ell \rightarrow \cdots \rightarrow k$ that d-connects $\ell$ and $k$ given $S \setminus \lbrace i, j, k \rbrace$ and $\ell \not\in S$.  
Such a path exists since $\ell$ is an ancestor of $k$ and not an ancestor of all other nodes in $S \setminus \lbrace i, j, k \rbrace$.    So in both cases $i$ and $k$ are also d-connected given $S \setminus \lbrace i, j, k \rbrace$ using a path that does not containing the node $j$. Hence, $i$ and $k$ are also d-connected given $S \setminus \lbrace i, k \rbrace$, a contradiction. 

Next, we prove for $ i, j \in \adj(G_S, k)$, if $(i,j)$ is not an edge in $G_S$, then $(i,j)$ is an edge in $G_{S \setminus \{k \}}$.
Since $i \in \adj(G_S, k)$, there exists a path $P_i$ that d-connects $i$ and $k$ given $S \setminus \{i, k\}$, and similar for $j$. 
Using the same argument as the above, $i$ and $j$ are also d-connected to $k$ using $P_i$ and $P_j$, respectively, given $S \setminus \{i, j, k\}$. Defining $P$ as the path that combines $P_i$ and $P_j$, then $k$ must be a {non-collider} along $P$ as otherwise $i$ and $j$ would be d-connected given $S \setminus \{i, j\}$, in which case $i$ and $j$ would also be d-connected given $S \setminus \{i, j, k\}$, and $(i,j)$ would be an edge in $G_{S \setminus \{k \}}$.
\end{proof}

\subsection{Proof of Theorem~\ref{thm: NBMDA is equivalent to MDA in oracle setting}}
\label{subsubsec: proof of Theorem MDA equivalence}
In the oracle setting, there are two main differences between Algorithm~\ref{alg:emmd} and the minimum degree algorithm.  
First, Algorithm~\ref{alg:emmd} uses partial correlation testing to construct a graph, while the minimum degree algorithm uses the precision matrix $\Theta$.  
The second difference is that Algorithm~\ref{alg:emmd} only updates based on the partial correlations of neighbors of the tested nodes.

Let $\Theta_S$ denote the precision matrix of the marginal distribution over the variables $\lbrace X_i : i \in S \rbrace$.  
Since the marginal distribution is Gaussian, the $(i,j)$-th entry of $\Theta_S$ is nonzero if and only if $\rho_{i,j \vert S \setminus \lbrace i,j \rbrace} \neq 0$.  
Thus, to prove that Algorithm~\ref{alg:emmd} and the minimum degree algorithm are equivalent, it suffices to show the following:
Let $G_S$ be an undirected graph with edges corresponding to the nonzero entries of $\Theta_S$. Then for any node $k$, the graph $G_{S \setminus \lbrace k \rbrace}$ constructed as defined in Algorithm~\ref{alg:emmd} has edges corresponding to the nonzero entries of $\Theta_{S \setminus \lbrace k \rbrace}$.
To prove that this is indeed the case, note that by Lemma~\ref{lm:mmd}, if $G_S$ is already estimated then nodes $i$ and $j$ are connected in $G_{S \setminus \lbrace k \rbrace}$ if and only if $\rho_{i,j \vert S \setminus \lbrace i, j, k \rbrace} \neq 0$.  
Finally, since the marginal distribution over $S$ is multivariate Gaussian, the $(i,j)$-th entry of $\Theta_{S \setminus \lbrace k \rbrace}$ is non-zero if and only if $\rho_{i,j \vert S \setminus \lbrace i, j, k \rbrace} \neq 0$.
\hfill$\square$

\subsection{Proof of Theorem~\ref{thm: NBMDA consistent in non-oracle setting}}
Let $\PP^{\textrm{oracle}}(\hat\pi)$ denote the probability that $\hat\pi$ is output by Algorithm~\ref{alg:emmd} in the oracle-setting, and let $N_{\hat\pi}$ denote the number of partial correlation tests that had to be performed.  
Then $N_{\hat\pi}\leq \mathcal{O}(p d_{\hat\pi}^2)$, where $d_{\hat\pi}$ is the maximum degree of the corresponding minimal independence map $\mathcal{G}_{\hat\pi}$.  
Therefore, using the same arguments as in the proof of Theorem~\ref{thm_high_dim_consist}, we obtain:
\begin{align*}
\begin{split}
\PP&[\hat\pi\; \text{is generated by Algorithm~\ref{alg:emmd}}]\\
& \geq \PP^{\textrm{oracle}}(\hat\pi) \PP[\text{all hypothesis tests for generating $\hat\pi$ are consistent}] \\
& \geq \PP^{\textrm{oracle}}(\hat\pi) \left(1 - \mathcal{O}(p d_{\hat\pi}^2) \underset{(i,j,S) \in N_{\hat\pi}}{\sup} \PP(E_{i,j \vert S}) \right), \\
& \geq \PP^{\textrm{oracle}}(\hat\pi) \left( 1 - \mathcal{O}(\exp(2 \log d_{\hat\pi} + \log p + \log n -c' n^{1 - 2\ell})) \right), \\
& \geq \PP^{\textrm{oracle}}(\hat\pi) \left( 1 - \mathcal{O}(\exp(-c n^{1 - 2\ell})) \right). \\
\end{split}
\end{align*}
Let $\Pi$ denote the set of all possible output permutations of the minimum degree algorithm applied to $\Theta$. Then
\begin{align*}
\begin{split}
\PP&[\text{Algorithm~\ref{alg:emmd} outputs a permutation in}\; \Pi] \\
& \geq \underset{\hat\pi \in \Pi}{\sum} \PP[\hat\pi\; \text{is output by Algorithm~\ref{alg:emmd}}], \\
& \geq 1 - \mathcal{O}(\exp(-c n^{1 - 2\ell})),
\end{split}
\end{align*}
which completes the proof.
\hfill$\square$

\begin{figure}[t!]
\centering
\includegraphics[width=0.4\textwidth]{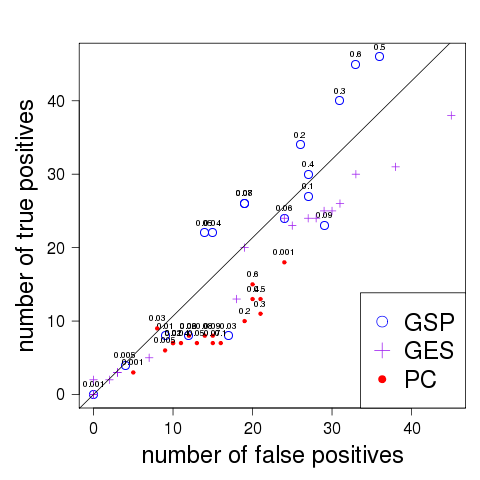}
\caption{Performance of the causal network learned by Algorithm~\ref{alg: greedy sp depth and start control} with $d = 4$ and $r = 20$ as compared to the PC-algorithm and GES in predicting the effect of each intervention using a q-value cutoff of $1$; line corresponds to random guessing. }
	\label{fig: Dixit supp}
	\vspace{-0.2cm}
\end{figure}

\section{Additional figures for experiments}

In this section, we present an additional figure supporting our experimental findings in Section~\ref{sec_real data}. Figure~\ref{fig: Dixit supp} shows the resulting receiver operating characteristic curves for the greedy sparsest permutation algorithm, the PC-algorithm as well as greedy equivalence search when using a q-value of 1 to identify true positive / false positive edges. More specifically, we consider an arrow from gene $A$ to gene $B$ in the learned network as a true positive if the magnitude of the corresponding q-value is larger than $1$, and a false positive otherwise. The random guessing line was adjusted accordingly. Our greedy sparsest permutation algorithm outperforms the PC-algorithm and greedy equivalence search, which both perform similar to random guessing.

\section{Computational Times for Simulations}

To test the computational efficiency of our greedy sparsest permutation algorithm, we compared its run time to the PC-algorithm and greedy equivalence search in the setting $p = 8, s = 4$ and $n = 1000$, which is the setting considered in Figures~\ref{fig: low dimensional ROCs} (a)-(b) in the main paper. For a fair comparison, selected the hyperparameters of each algorithm so that the resulting graphs have a similar sparsity, namely 0.001 for our greedy sparsest permutation search, 0.01 for the PC-algorithm and $\lambda_n=1/2 \log(n)$ for greedy equivalence search. 
The R implementation of our greedy sparsest permutation algorithm used in this paper took $0.42$ seconds for one run, while it took  $0.08$ seconds for the  PC-algorithm and $0.02$ seconds for greedy equivalence search (using the pcalg package in R). While our implementation of the greedy sparsest permutation algorithm should be seen mainly as a proof-of-concept, in the meantime, a faster implementation of the greedy sparsest permutation algorithm has been developed and is available as a python package at \url{https://github.com/uhlerlab/causaldag}.

Finally, we note that the moves used in Algorithm~\ref{alg: greedy sp depth and start control} are a strict subset of the moves used by the algorithm of \cite{TK12}.
Moreover, this subset explicitly excludes moves that are guaranteed not to improve the value of the score function.  
Therefore, it seems likely that  Algorithm~\ref{alg: greedy sp depth and start control} performs with efficiency comparable or favorable to the algorithm of \cite{TK12}, which was already shown to be more efficient than the greedy equivalence search.

\end{appendix}

\end{document}